\documentclass[12pt]{article}
\usepackage{latexsym}
\usepackage{epsfig,amssymb,euscript}
\usepackage{amsmath}
\numberwithin{equation}{section}  

\newsavebox{\ns}
\newsavebox{\dbrane}
\newsavebox{\dbshort}

\def\be{\begin{equation}}
\def\ee{\end{equation}}
\def\bea{\begin{eqnarray}}
\def\eea{\end{eqnarray}}

\newcommand{\nn}{\nonumber}

\def\Dslash{\,\,{\raise.15ex\hbox{/}\mkern-12mu D}}
\def\Dbarslash{\,\,{\raise.15ex\hbox{/}\mkern-12mu {\bar D}}}
\def\delslash{\,\,{\raise.15ex\hbox{/}\mkern-9mu \partial}}
\def\delbarslash{\,\,{\raise.15ex\hbox{/}\mkern-9mu {\bar\partial}}}
\def\pslash{\,\,{\raise.15ex\hbox{/}\mkern-9mu p}}
\def\calDslash{\,\,{\raise.15ex\hbox{/}\mkern-12mu {\cal D}}}

\newcommand\numax{\nu_{\mathrm{max}}}

\newcommand\R{\mathbb{R}}
\newcommand\Z{\mathbb{Z}}
\newcommand\F{\mathbb{F}}

\newcommand\C{\mathbb{C}}
\newcommand\T{\mathbb{T}}
\newcommand\diff{\mathrm{d}}

\newtheorem{theorem}{Theorem}[section] 
\newtheorem{lemma}[theorem]{Lemma} 
\newtheorem{proposition}[theorem]{Proposition} 
\newtheorem{corollary}[theorem]{Corollary}

\newenvironment{proof}[1][Proof]{\begin{trivlist} 
\item[\hskip \labelsep {\bfseries #1}]}{\end{trivlist}}

\newcommand{\qed}{\nobreak \ifvmode \relax \else \ifdim\lastskip<1.5em \hskip-\lastskip \hskip1.5em plus0em minus0.5em \fi \nobreak \vrule height0.75em width0.5em depth0.25em\fi}

\begin{document}
\begin{titlepage}
\begin{center}

\vspace* {2.5 cm} 
{\Large \bf Resolutions of non-regular}\\
\vskip 4mm
{\Large \bf Ricci-flat K\"ahler cones}\\

\vskip 1cm
{Dario Martelli$^{1*}$ and James Sparks$^{2,3}$}\\
\vskip 1cm

1: \emph{Institute for Advanced Study\\
Einstein Drive, Princeton, NJ 08540, U.S.A.}\\
\vskip 0.5cm
2: \emph{Department of Mathematics, Harvard University \\
One Oxford Street, Cambridge, MA 02138, U.S.A.}\\
\vskip 0.5cm
3: \emph{Jefferson Physical Laboratory, Harvard University \\
Cambridge, MA 02138, U.S.A.}\\

\vskip 1cm

\end{center}

\begin{abstract}
\noindent  We present explicit constructions of complete Ricci-flat 
K\"ahler metrics that are asymptotic to cones over non-regular 
Sasaki-Einstein manifolds. The metrics are constructed from a complete 
K\"ahler-Einstein manifold $(V,g_V)$ of 
positive Ricci curvature and admit a Hamiltonian two-form of order 
two. 
We obtain Ricci-flat K\"ahler metrics on the total spaces of (i) 
holomorphic $\C^2/\Z_p$ orbifold fibrations over $V$, 
(ii) holomorphic orbifold fibrations over weighted projective 
spaces $\mathbb{WCP}^1$, with generic fibres being the canonical complex cone over $V$, 
and (iii) the canonical orbifold line bundle over a family of Fano orbifolds.
As special cases, we also obtain smooth complete Ricci-flat K\"ahler metrics 
on the total spaces of (a) rank two holomorphic vector bundles over $V$, and  
(b) the canonical line bundle over a family of geometrically ruled Fano manifolds 
with base $V$. 
When $V=\mathbb{CP}^1$ our results give Ricci-flat K\"ahler orbifold metrics on 
various toric partial resolutions of the cone over the Sasaki-Einstein manifolds $Y^{p,q}$.
\end{abstract}

\vfill 
\hrule width 5cm
\vskip 5mm

{\noindent $^*$ {\small On leave from: \emph{Blackett Laboratory, 
Imperial College, London SW7 2AZ, U.K.}}}

\end{titlepage}


\pagestyle{plain}
\setcounter{page}{1}
\newcounter{bean}
\baselineskip18pt

\tableofcontents

\section{Introduction and summary}

\subsection{Introduction}

A Sasaki-Einstein manifold $(L,g_L)$ is a complete Riemannian manifold 
whose metric cone
\bea\label{cone}
C(L)\, =\, \R_+\times L~, \qquad  g_{C(L)} \, =\,  \diff r^2 + r^2 g_L\eea
is Ricci-flat K\"ahler. 
The metric in (\ref{cone}) 
is singular at $r=0$, unless $(L,g_L)$ is the round sphere, and it 
is natural to ask whether there exists a resolution 
{\it i.e.} a complete Ricci-flat K\"ahler metric on a non-compact manifold 
$X$ which is asymptotic to the cone (\ref{cone}). 
More generally, one 
can consider \emph{partial} resolutions in which $X$ also has 
singularities. There are particularly strong physical motivations for studying 
such partial resolutions; for example, certain types of orbifold 
singularity are well-studied in String Theory, and may give rise to 
interesting phenomena, such as non-abelian gauge symmetry.

These geometrical structures are of particular interest in the AdS/CFT 
correspondence \cite{Maldacena}. In complex dimension three 
or four, a Ricci-flat K\"ahler cone $C(L)$ is AdS/CFT dual to a supersymmetric conformal field theory
in dimension four or three, respectively. 
Resolutions of such conical singularities are then of interest for a number of 
different physical applications. For example, in AdS/CFT such resolutions 
correspond to certain deformations of the conformal field theory.

On a K\"ahler cone $(C(L),g_{C(L)})$ there is a canonically defined vector field, the 
Reeb vector field:
\be
\xi = J\left(r\frac{\partial}{\partial r}\right)
\ee
where $J$ denotes the complex structure tensor on the cone. $\xi$ 
is a holomorphic Killing vector field, and has unit norm on the link 
$L=\{r=1\}$ of the singularity at $r=0$. If the orbits 
of $\xi$ all close then $\xi$ generates a $U(1)$ isometry of $(L,g_L)$, which 
necessarily acts locally freely since $\xi$ is nowhere zero, 
and the Sasakian structure is said to be 
either regular or quasi-regular if this action is free or not, respectively. 
The orbit space is in general a 
K\"ahler-Einstein orbifold $(M,g_M)$ of positive Ricci curvature, which is a smooth 
manifold in the regular case. More generally, 
the orbits of $\xi$ need not all close, in which case 
the Sasakian structure is said to be irregular. 

Suppose that $(L,g_L)$ is a \emph{regular} Sasaki-Einstein manifold. 
In this case $L$ is a $U(1)$ fibration over a K\"ahler-Einstein 
manifold
$(M,g_M)$, which we assume\footnote{$b_1(M)=0$ necessarily 
 \cite{kobayashi}.} is 
simply-connected. 
Let $K_M$ denote the canonical line bundle of 
$M$, and let $I$ denote the Fano index of $M$. The latter is the largest 
positive integer such that
\bea
K_M^{-1/I}\in \mathrm{Pic}(M) = H^2(M;\Z)\cap H^{1,1}(M;\C)~.\eea
It is then well-known that the simply-connected 
cover of $L$ is diffeomorphic to the unit circle bundle in the holomorphic 
line bundle
$K_M^{1/I}$. Taking the quotient of $L$ by $\Z_m\subset U(1)$ gives 
instead a smooth Sasaki-Einstein manifold diffeomorphic to the 
unit circle bundle in $K_M^{m/I}$. For example, suppose $(M,g_M)$ is 
$\mathbb{CP}^2$ 
equipped with its Fubini-Study metric. Then the Fano index is $I=3$, 
and the canonical line 
bundle is $K_{\mathbb{CP}^2}={\cal O}(-3)$. The total space of 
the associated circle bundle is thus $S^5/\Z_3$, whereas
the simply-connected cover of  $(L,g_L)$ is $S^5$ equipped with its round metric.

When $m=I$, there is a canonical way of resolving the above 
Ricci-flat K\"ahler cone: there exists a smooth complete 
Ricci-flat K\"ahler metric on the total 
space of the canonical line bundle $K_M$ over $M$. 
The metric is in fact explicit, up to the 
K\"ahler-Einstein metric $g_M$ on $M$, and is constructed using the 
Calabi ansatz \cite{Calabi}. These metrics were constructed in the mathematics literature in 
\cite{BB}, and in the 
physics literature in \cite{PP}. More generally, 
there may exist other resolutions. The simplest example is perhaps given by 
$M=\mathbb{CP}^1\times\mathbb{CP}^1$ (also known as zeroth Hirzebruch 
surface, denoted 
$\mathbb{F}_0$) with its standard 
K\"ahler-Einstein metric. Here $I=2$, and for $m=2$ the construction of \cite{BB,PP} produces 
a complete metric on the total space of $K_{\F_0}$, which is 
asymptotic to a cone 
over the homogeneous Sasaki-Einstein manifold $T^{1,1}/\Z_2$. On 
the other hand, the cone over the Sasaki-Einstein manifold with 
$m=1$ instead has a small resolution: there is a 
smooth complete Ricci-flat K\"ahler metric 
on the total space of the 
rank two holomorphic vector bundle $\mathcal{O}(-1)\oplus \mathcal{O}(-1)$ 
over $\mathbb{CP}^1$, which is asymptotic to a cone over $T^{1,1}$. This is known in the physics literature as 
the resolved conifold metric \cite{candelas}.

More generally, there are the existence results of Tian and Yau 
\cite{tianyau1, tianyau2}. In the latter reference it is proven that, 
under certain mild assumptions, $X=\bar{X}\setminus D$ admits 
a complete Ricci-flat K\"ahler metric that is asymptotic to a 
cone, provided that the divisor $D\subset \bar{X}$ in the compact 
K\"ahler manifold (or orbifold) $\bar{X}$ admits a K\"ahler-Einstein 
metric of positive Ricci curvature. 
These metrics are therefore also asymptotic to cones over regular, 
or quasi-regular, Sasaki-Einstein manifolds. However, the metrics that 
we shall present  in this paper lie outside this class, and 
their existence was not guaranteed by any theorem.

In \cite{paper1, paper2, paper3} infinite families 
of explicit Sasaki-Einstein manifolds were constructed, in all odd 
dimensions, in both 
the quasi-regular and irregular classes. In particular, these 
were the first examples of irregular Sasaki-Einstein manifolds. 
The construction produces, for each complete K\"ahler-Einstein 
manifold $(V,g_V)$ of positive Ricci curvature, an infinite family $Y^{p,k}(V)$ of associated 
complete Sasaki-Einstein manifolds. Here $p$ and $k$ are positive 
integers satisfying $pI/2<k<pI$, where $I$ is the Fano index of $V$.
Given the above results, it 
is natural to investigate whether or not there exist resolutions
of the corresponding 
Ricci-flat K\"ahler cones. In fact examples of 
such resolutions have recently been 
constructed in \cite{japanese} and \cite{lupope}. In this paper we 
significantly generalise these results; 
the results of \cite{japanese, lupope}
are recovered by substituting $(V,g_V)=\mathbb{CP}^1$ (in particular, $I=2,m=1$) with its standard 
metric, or $(V,g_V)=$ product 
of complex projective spaces, into Corollary \ref{ruled}, respectively.

Quite recently, Futaki \cite{futtaki} has used the Calabi ansatz 
to construct complete
Ricci-flat K\"ahler metrics on the canonical line bundles ({\it i.e.} 
$m=I$, in the above notation) over toric Fano manifolds. 
A key point in the construction is the general existence result 
of \cite{toricSE} for toric Sasaki-Einstein metrics on links 
of isolated toric Gorenstein singularities.

\subsection{Summary}

Our constructions are based on a class of explicit
local K\"ahler metrics that have appeared recently in the mathematics 
literature \cite{apostolov1, apostolov2, apostolov3}
and have been independently discovered in the physics literature in 
\cite{chen}.
The metrics we study all admit a \emph{Hamiltonian two-form}, in the sense of \cite{apostolov1}, of order two.  As noted in 
\cite{apostolov1}, the Calabi ansatz is a special case of a local 
K\"ahler metric admitting a Hamiltonian two-form of order one. 
More generally, a K\"ahler metric admitting a 
Hamiltonian two-form of order one locally fibres over a product 
of K\"ahler manifolds: this ansatz was in fact used in 
the paper \cite{gibbons} to construct complete Ricci-flat K\"ahler 
metrics on various holomorphic vector bundles over products of 
K\"ahler-Einstein manifolds; the asymptotic cones are again all 
regular, however. For simplicity, we study here only a single 
K\"ahler-Einstein manifold $(V,g_V)$, rather than a product of 
K\"ahler-Einstein manifolds\footnote{The product case was 
in fact discussed briefly in \cite{paper3}, with some global 
analysis of the corresponding Sasaki-Einstein metrics appearing in 
\cite{otherpope}.}.
The local metrics depend on two real parameters. In sections \ref{asymptotics} and 
\ref{global}  we establish that
it is possible to choose one metric parameter $\nu$ in such a way that the 
metric asymptotes to a cone over one of the
non-regular Sasaki-Einstein metrics constructed in \cite{paper3}; there 
are a countably infinite number of such choices for $\nu$. 
The remainder of
the paper is then devoted to analysing in detail the various 
possibile choices for the second metric parameter 
$\mu$.  We obtain
Ricci-flat K\"ahler metrics on \emph{partial resolutions}, with  
singularities that we carefully describe, as well as various
\emph{smooth complete Ricci-flat K\"ahler} metrics, that provide
 distinct resolutions of the conical singularities. 
When the Fano manifold $(V,g_V)$
 is toric, the resolutions we construct are all toric resolutions. 
In particular, when $V=\mathbb{CP}^1$ equipped with its 
standard round metric, our results may be described in terms of 
the toric geometry of the family $C(Y^{p,q})$ of isolated toric Gorenstein 
singularities \cite{toricpaper}. 
Such a description, together with the AdS/CFT 
interpretation of the metrics constructed here, will appear 
elsewhere \cite{np2}.

In section \ref{smallsection} we investigate two classes 
of (partial) resolution 
that we shall refer to as \emph{small resolutions}.
This nomenclature is motivated by the fact that these metrics may 
be thought of as two different generalisations
of the resolved  conifold metric on 
${\cal O}(-1)\oplus {\cal O}(-1)\to \mathbb{CP}^1$. 
First, we obtain complete 
Ricci-flat K\"ahler orbifold metrics on the total spaces of certain 
holomorphic $\C^2/\Z_p$ fibrations over $V$. When $p=1$ this leads to
smooth  Ricci-flat K\"ahler metrics on the total spaces of certain 
rank two holomorphic vector bundles over the Fano $V$, as summarised 
in Corollary \ref{p=1} below.
For instance, taking $V=\mathbb{CP}^2$ with its standard 
K\"ahler-Einstein metric, we obtain a smooth
complete metric on the total space of the rank two holomorphic 
vector bundle ${\cal O}(-2)\oplus {\cal O}(-1)\to \mathbb{CP}^2$.
On the other hand, these results also produce an infinite family 
of \emph{partial} small resolutions of the cones 
over the non-regular 
Sasaki-Einstein manifolds $Y^{p,q}$ \cite{paper2}. The resolution is in general 
only partial since the normal fibre to the blown-up 
$\mathbb{CP}^1$ is $\C^2/\Z_p$. The fibers are 
non-trivially twisted over 
$\mathbb{CP}^1$, with the form of the twisting depending on the integer $q$. 
When $p=1$ we recover precisely the resolved conifold metric. 
In section \ref{floppy} we will describe a second type of partial 
resolution, 
where one instead blows up a weighted projective space. 
The partial resolution is a fibration over this weighted projective space 
with generic fibres being the singular 
canonical complex cone over the Fano $V$ ({\it cf}. Theorem \ref{flopper} below).
 In particular, when $V=\mathbb{CP}^1$, the latter fibres
are simply copies of $\C^2/\Z_2$, which is the blow-down of 
${\cal O}(-2)\to \mathbb{CP}^1$. More precisely, in this case we obtain a 
$\C^2/\Z_2$ orbifold fibration over $\mathbb{WCP}^1_{[d,p-d]}$, 
where $d=k/2$ implies that $k=p+q$ must be even.

In section \ref{bigones} we investigate a class of complete 
Ricci-flat K\"ahler orbifold metrics on the total space of the 
canonical line bundle over a family of Fano orbifolds. 
These are a direct generalisation of the work of \cite{BB,PP}, 
which was based on the Calabi ansatz. 
Our Fano orbifolds are $\mathbb{WCP}^1_{[r,p-r]}$ fibrations over $V$, 
where $0<r<k/I$. 
The induced orbifold metric on $M$, which is the zero-section of 
the canonical line bundle, is K\"ahler, but $(M,g_M)$ is \emph{not} 
K\"ahler-Einstein. $M$ is smooth if and only if $p=2$, $r=1$ and in this 
case $M$ is a $\mathbb{CP}^1$ fibration over the Fano $V$, of the form 
$M=\mathbb{P}_V(\mathcal{O}\oplus K_V^{m/I})$ where $0<m<I$. 
For instance, when $V=\mathbb{CP}^1$, $M$ is the first del
Pezzo surface, which is well known to have non-vanishing Futaki invariant, 
the latter being 
an obstruction to the existence of a K\"ahler-Einstein metric. 
More generally we obtain smooth complete metrics
on the total space of the canonical line bundle $K_M$ over $M$, 
generalising \cite{BB,PP}
to the case of non-regular Sasaki-Einstein boundaries.


We summarise our results more formally by the following Theorems. 
Note that for general $p$ and $k$ the Sasaki-Einstein manifolds 
 $Y^{p,k}(V)$ \cite{paper3} are irregular:

\begin{theorem} \label{orbifolds} Let $(V,g_V)$ be a complete 
K\"ahler-Einstein manifold of positive 
Ricci curvature with canonical line bundle $K_V$ and Fano index $I$. 
Then for every $p, k\in \mathbb{N}$
positive integers with $pI/2 < k < pI$ there is an explicit 
complete Ricci-flat K\"ahler orbifold metric on the total space of a $\C^2/\Z_p$ bundle
over $V$. Here $\Z_p\subset U(1)\subset SU(2)$ acts on $\C^2$ in the standard 
way, and the bundle is given by

\bea
\left[K_V\oplus K_V^{k/I}\right]\times_{\lambda} 
\C^2/\Z_p
\eea
where 
\bea\label{action}
\lambda :  S^1 \times S^1 \times \C^2/\Z_p  & \rightarrow & \C^2/\Z_p\nonumber\\
(\theta_1,\theta_2;z_1,z_2) & \mapsto & (\exp(i\theta_1-i\theta_2/p)z_1, 
\exp(i\theta_2/p)z_2)
\eea
and $z_1,z_2$ are standard complex coordinates on $\C^2$. 
The metric asymptotes to a cone over the Sasaki-Einstein
manifold $Y^{p,k}(V)$. \end{theorem}

When $p=1$ we obtain a finite number of completely smooth resolutions, for 
each $(V,g_V)$. These may be regarded as higher-dimensional versions 
of the small resolution of the conifold, which are asymptotic to 
non-regular Ricci-flat K\"ahler cones. Setting $p=1$, $m=I-k$ in Theorem \ref{orbifolds} gives

\begin{corollary}
\label{p=1} Let $(V,g_V)$ be a complete K\"ahler-Einstein manifold of positive 
Ricci curvature
with canonical line bundle $K_V$ and Fano index $I$. Then for every 
$m\in\mathbb{N}$ with
$0<m<I/2$ there is an explicit smooth complete asymptotically conical
Ricci-flat K\"ahler metric on the total space of the rank two holomorphic
vector bundle $K_V^{m/I}\oplus K_V^{(I-m)/I}$ over $V$. 
The metric asymptotes to a cone
over the Sasaki-Einstein manifold $Y^{1,I-m}(V)$. 
\end{corollary} 

We also obtain 

\begin{theorem}
\label{flopper}
Let $(V,g_V)$ be a complete K\"ahler-Einstein
manifold of positive Ricci curvature with canonical line bundle $K_V$ 
and Fano index $I$. Then for each $p,d\in\mathbb{N}$ with $p/2<d<p$ 
there is an explicit complete Ricci-flat K\"ahler orbifold metric on the 
total space of the canonical 
complex cone $\C_V$ over $V$, fibred over the weighted 
projective space $\mathbb{WCP}^1_{[d,p-d]}$. The fibration structure 
is given by the orbifold fibration
\bea
K_{\mathbb{WCP}^1_{[d,p-d]}}\times_{U(1)} \C_V~.\eea
Here the $U(1)\subset\C^*$ action is the standard one on the complex cone 
$\C_V$. The metric is 
completely smooth away from the tip of the complex cone fibres, and asymptotes
to a cone over the Sasaki-Einstein manifold $Y^{p,Id}(V)$.
\end{theorem}

In section \ref{bigones} we prove

\begin{theorem} \label{bigorbifolds} Let $(V,g_V)$ be a complete 
K\"ahler-Einstein manifold of positive 
Ricci curvature with canonical line bundle $K_V$ and Fano index $I$. 
Then for every $p, k, r\in \mathbb{N}$
positive integers with $p/2 < k/I < p$, $0<r<k/I$, 
there is an explicit smooth complete
Ricci-flat K\"ahler orbifold metric on the total space of the 
canonical line bundle $K_M$ over the Fano orbifold
\bea
M = K_V^{m/I}\times_{U(1)}\mathbb{WCP}^1_{[r,p-r]} ~,\eea
where $m=k-rI$. Here we use the standard effective action of $U(1)$ on the 
weighted projective space $\mathbb{WCP}^1_{[r,p-r]}$, with 
orientation fixed so that the section with normal fibre $\C/\Z_{p-r}$ 
has normal bundle $K_V^{m/I}$.
The metric asymptotes, for every $r$, to a cone over the Sasaki-Einstein
manifold $Y^{p,k}(V)$. 
\end{theorem}

Setting $p=2$, $r=1$ in Theorem \ref{bigorbifolds} effectively blows 
up the zero section of the orbifold metric in Theorem 
\ref{orbifolds} to again obtain a finite number of completely smooth 
resolutions, for each $(V,g_V)$:

\begin{corollary} 
\label{ruled} 
Let $(V,g_V)$ be a complete K\"ahler-Einstein manifold of 
positive Ricci curvature
with canonical line bundle $K_V$ and Fano index $I$. Then for each 
$m\in\mathbb{N}$ with $0<m<I$
there is an explicit smooth complete  
Ricci-flat K\"ahler metric on the total space of the canonical line bundle $K_M$ over the geometrically
ruled Fano manifold $M = \mathbb{P}_V(\mathcal{O}\oplus K_V^{m/I})$. 
The metric asymptotes to a cone
over the Sasaki-Einstein manifold $Y^{2,m+I}(V)$. \end{corollary}

We note that $I\leq n+1$ with equality if and only if $V=\mathbb{CP}^n$. 
In fact also $I=n$ if and only if $V=\mathbb{Q}^n$ is the quadric 
in $\mathbb{CP}^{n+1}$ -- see, for example, \cite{Kollarbook}. Both of these examples 
admit homogeneous K\"ahler-Einstein metrics.


\section{Local metrics}

In this section we introduce the class of explicit 
local K\"ahler metrics that we wish to study. 
These metrics all admit a \emph{Hamiltonian two-form}
\cite{apostolov1}. In section \ref{localmet} we give a brief review of local 
K\"ahler metrics admitting Hamiltonian two-forms, focusing on the relevant 
cases of order one and order two, and present the local form of the metrics 
used throughout the remainder of the paper. In section \ref{complex} 
we introduce local complex coordinates. Finally, section \ref{asymptotics} 
demonstrates that, in a certain limit, the local metrics are 
asymptotic to a cone over the local class of Sasaki-Einstein metrics 
studied in \cite{paper3}.

\subsection{K\"ahler metrics with Hamiltonian two-forms}\label{localmet}

If $(X,g,J,\omega)$ is a K\"ahler structure, then a \emph{Hamiltonian 
two-form} $\phi$ is a real $(1,1)$-form that solves non-trivially the equation 
\cite{apostolov1}
\bea
\nabla_Y \phi \, = \, \frac{1}{2}\left(\diff \ \mathrm{tr}_{\omega}\phi 
\wedge JY^{\flat} - J\diff \ \mathrm{tr}_{\omega}\phi\wedge Y^{\flat}\right)~.
\eea
Here $Y$ is any vector field, $\nabla$ denotes the Levi-Civita connection, 
and $Y^{\flat}=g(Y,\cdot)$ is the one-form dual to $Y$. 

The key result for our purposes is that 
the existence of $\phi$ leads to an ansatz for the K\"ahler metric $g$ 
such that Ricci-flatness is equivalent to solving a simple set of decoupled
ordinary differential equations. We therefore merely sketch 
the basic ideas that lead to this result; for a full exposition 
on Hamiltonian two-forms, the reader is referred to \cite{apostolov1}. 
We note that many of these ans\"atze had been 
arrived at prior to the work of \cite{apostolov1}, both 
in the mathematics literature (as pointed out in \cite{apostolov1}), 
and also in the physics literature. The theory of Hamiltonian two-forms 
unifies these various approaches.

One first notes that if $\phi$ is a Hamiltonian two-form, then 
so is $\phi_t = \phi - t\omega$ for any $t\in\mathbb{R}$. One then 
defines the \emph{momentum polynomial} of $\phi$ to be
\bea
p(t) \, = \, \frac{(-1)^N}{N!} *\phi_t^N~.
\eea
Here $N$ is the complex dimension of the K\"ahler manifold and $*$ is 
the Hodge operator with respect to the metric $g$. It is then 
straightforward to show that $\{p(t)\}$ are a set of 
Poisson-commuting Hamiltonian functions for the one-parameter family of 
Killing vector fields $K(t)=J\mathrm{grad}_g p(t)$. For a fixed point 
in the K\"ahler manifold, these Killing vectors will span a vector 
subspace of the tangent space of the point; the maximum dimension of 
this subspace, taken over all points, is called the \emph{order} $s$ of 
$\phi$. This leads to a Hamiltonian $\mathbb{T}^s$ action, at least locally, 
on the K\"ahler 
manifold, and one may then take a (local) K\"ahler quotient by 
this torus action. The reduced K\"ahler metric depends on the 
moment map level at which one reduces, but only very weakly: 
the reduced K\"ahler metric is a direct product of $S$ K\"ahler 
manifolds $(V_i,c_i(\mathbf{\mu})g_{V_i})$, $i=1,\ldots, S$, where $c_i(\mathbf{\mu})$ are 
functions of the moment map coordinates $\mathbf{\mu}$. The $2s$-dimensional 
fibres 
turn out to be \emph{orthotoric}, which is a rather special type  
of toric K\"ahler structure. For further details, we refer the reader to 
reference \cite{apostolov1}. 

The simplest non-trivial case is a Hamiltonian two-form of order one, 
with $S=1$. This turns out to be precisely the \emph{Calabi ansatz} 
\cite{Calabi}. The local metric and K\"ahler form may be written in the form
\bea\label{calabiform}
g & = & (\beta-y)g_V + \frac{\diff y^2}{4Y(y)}+Y(y)(\diff \psi + A)^2~,\nonumber\\
\omega & = & (\beta-y)\omega_V - \frac{1}{2}\diff y\wedge (\diff\psi+A)~.
\eea
Here $A$ is a local 1-form on the K\"ahler manifold $(V,g_V, \omega_V)$ 
satisfying $\diff A = 2\omega_V$. The Killing vector field $\partial/\partial\psi$ 
generates the Hamiltonian action, with Hamiltonian function $y$. 
The momentum polynomial is given by $p(t)=(t-y)(t-\beta)^{N-1}$ 
where $\beta\in\R$ is a constant. 
Calabi used this ansatz to produce an explicit family of so-called extremal 
K\"ahler metrics on the blow-up of $\mathbb{CP}^2$ at a point. One 
of these metrics is conformal to Page's Einstein metric \cite{Page}, 
which is perhaps more well-known to physicists. The same ansatz 
was used in \cite{BB, PP} to produce explicit constructions of 
complete non-compact K\"ahler metrics; indeed, this leads to the
construction of complete Ricci-flat K\"ahler metrics on $K_M$, 
where $(M,g_M)$ is a complete K\"ahler-Einstein manifold of positive 
Ricci curvature. The general form of a K\"ahler 
metric with a Hamiltonian two-form of order one allows one to 
replace $(V,g_V)$ by a direct product of $S>1$ K\"ahler manifolds, as 
mentioned above. In fact precisely this ansatz was used in 
section 5 of \cite{gibbons}, before the work of 
\cite{apostolov1}, to produce a number of examples of 
complete non-compact Ricci-flat K\"ahler manifolds. The same general 
form 
was thoroughly investigated in \cite{apostolov3}, and used to 
give explicit constructions of compact extremal K\"ahler manifolds.

In this paper we study the case of a Hamiltonian two-form of order two, 
with $S=1$. A K\"ahler structure $(X,g,\omega)$ admitting such a two-form may be written 
in the form
\bea\label{metric}
g & = & \frac{(\beta-x)(\beta-y)}{\beta}g_V + \frac{y-x}{4X(x)}\diff x^2 + \frac{y-x}{4Y(y)}\diff y^2 \nonumber\\ && + \frac{X(x)}{y-x}\left[\diff\tau + \frac{\beta-y}{\beta}(\diff\psi + A)\right]^2  + \frac{Y(y)}{y-x}\left[\diff\tau + \frac{\beta-x}{\beta}(\diff\psi + A)\right]^2~,\\
\omega & = & \frac{(\beta-x)(\beta-y)}{\beta}\omega_V - 
\frac{1}{2}\diff x\wedge\left[\diff\tau + \frac{\beta-y}{\beta}(\diff\psi + A)\right]\nonumber\\ && - \frac{1}{2}\diff y \wedge\left[\diff\tau + \frac{\beta-x}{\beta}(\diff\psi + A)\right]~.\eea
Here $(V,g_V,\omega_V)$ is again a K\"ahler manifold with, locally, 
$\diff A = 2\omega_V$. The momentum polynomial is now given by 
$p(t) = (t-x)(t-y)(t-\beta)^n$, where we denote $n=\dim_{\C}V=N-2$. 
The Hamiltonian action is generated by the Killing vector fields 
$\partial/\partial\tau$, $\partial/\partial\psi$.

A computation shows that the metric (\ref{metric}) is Ricci-flat  
if $(V,g_V)$ is a K\"ahler-Einstein manifold of positive\footnote{More 
generally one might also 
consider zero or negative Ricci curvature.} Ricci curvature 
and the metric functions are given by
\bea\label{functions}
X(x) & = & \beta(x-\beta)+\frac{n+1}{n+2}c(x-\beta)^2+\frac{2\mu}{(x-\beta)^n}\nonumber\\ 
Y(y) & = & \beta(\beta-y)-\frac{n+1}{n+2}c(\beta-y)^2-\frac{2\nu}{(\beta-y)^n}~.\eea
Here $\beta$, $c$, $\mu$ and $\nu$ are real constants and, without loss of generality, 
we have normalised the metric $g_V$ so that $\mathrm{Ric}_V=2(n+1)g_V$.

Note that, provided $\beta\neq 0$, one may define $x=\beta\hat{x}$, $y=\beta\hat{y}$, multiply $g$ by $1/\beta$, and then 
relabel $\hat{x}=x$, $\hat{y}=y$ to obtain (\ref{metric}) with $\beta=1$. 
Similarly, provided $c\neq 0$, one may define $x^{\prime} = 1+c(x-1)$, $y^{\prime}=1+c(y-1)$, $\tau^{\prime} = c\tau$, multiply $g$ by 
$c^2$, and then relabel $x^{\prime}=x$, $y^{\prime}=y$, $\tau^{\prime}=\tau$ to obtain 
(\ref{metric}) with $c=1$. The cases $c=0$ and $\beta=0$ (accompanied by 
a suitable scaling of
the coordinates) are treated in the appendix, where the parameter $\beta$
is also further discussed. Henceforth we set $\beta=c=1$.

\subsection{Complex structure}\label{complex}

In this section we introduce a set of local complex coordinates 
on the (local) K\"ahler manifold $(X,g,\omega)$.
We first define the complex one-forms
\bea
\eta_1 &  = & \frac{\diff x}{2X(x)} -  \frac{\diff y}{2Y(y)} -
i\diff \psi\nn\\
\eta_2 & = &  \frac{1-x}{2X(x)}\diff x -  \frac{1 - y}{2Y(y)} \diff
y + i\diff \tau~.
\eea
The following is then a closed $(n+2,0)$-form:
\bea
\Omega \, =\, \kappa\sqrt{X(x)Y(y)} \left[(1- x)(1- y )\right]^{n/2} \, (\eta_1 - iA)\wedge 
\eta_2 \wedge \Omega_V\eea
where
\bea 
\kappa \, =\, \exp [i(n+1)(\tau+\psi)]
\eea
and $\Omega_V$ is the $n$-form on $V$ satisfying
\bea
\diff \Omega_V \,=\, i(n+1)A\wedge \Omega_V~.
\label{errr}
\eea
More precisely, we may introduce local complex coordinates $z_1,\dots, z_n$ 
on $V$ and locally write
\bea
\Omega_V \, = \, f_V \diff z_1 \wedge\cdots \wedge\diff z_n~.
\eea
Globally, $f_V$ is a holomorphic section of the 
anti-canonical line bundle of $V$; on the overlaps of local 
complex coordinate patches this transforms oppositely to $\diff z_1 \wedge 
\dots \wedge \diff z_n$, giving a globally defined $n$-form $\Omega_V$ 
on $V$. So $f_V\in H^0(V,K^{-1}_V)$. The 
holomorphicity of $f_V$ may be seen by comparing with
(\ref{errr}), which implies
\bea\label{hol}
\left[\diff \log f_V - i(n+1)A\right]\wedge \Omega_V & = & 0~.
\eea 
The local one-form $(n+1)A$ is a connection on the holomorphic 
line bundle $K^{-1}_V$, since 
$(n+1)\diff A = \rho_V$ is the Ricci form of $V$. Equation (\ref{hol}) 
then says that $f_V$ is a holomorphic section. Note that 
$\Omega_V$ has constant norm, even though $f_V$ necessarily has 
zeroes on $V$.

We may then introduce the local complex coordinates
\bea
Z_1 &= & \exp{\left[-i\psi + \int \frac{\diff x}{2X(x)}  - \frac{\diff y}{2Y(y)}\right]} \, f_V^{-1/(n+1)}\nn\\
Z_2 & = & \exp{\left[i\tau + \int \frac{(1 - x)\diff x}{2X(x) }  - \frac{(1 - y )\diff y}{2Y(y)}\right]}
\label{localZ}
\eea
satisfying
\bea
\diff \log Z_1 \,=\, \eta_1 - \frac{1}{(n+1)}\diff \log f_V~,\qquad \quad \diff \log Z_2 \,=\, \eta_2 ~.
\eea

\subsection{Asymptotic structure}
\label{asymptotics}

The metric (\ref{metric}) is symmetric in $x$ and $y$. We shall later 
break this symmetry by choosing one coordinate
to be a radial coordinate and the other to be a polar coordinate. 
Without loss of generality, we may take $x$ to be the radial coordinate. 
We analyse the metric in the limit $x\rightarrow \pm\infty$. 
Setting 
\bea
x\,=\, \pm\frac{n+1}{n+2}r^2
\eea
we obtain
$\mp g\rightarrow \diff r^2 + r^2 g_L$ where $g_L$ is the Sasaki-Einstein metric 
\bea\label{SEmetric}
g_L \,=\, \left(\frac{n+1}{n+2}\diff\tau+\sigma\right)^2 + g_T~.\eea 
Note in particular that for $x\rightarrow \infty$ it is $-g$ that 
is positive definite. In (\ref{SEmetric}) we have defined
\bea
\sigma \,=\, \frac{n+1}{n+2}(1-y)(\diff\psi+A)\eea
with $\diff\sigma=2\omega_T$, and $g_T$ is a local K\"ahler-Einstein metric 
given by 
\bea\label{gtrans}
g_T \,=\, \frac{n+1}{n+2}\left[(1-y)g_V + \frac{\diff y^2}{4Y(y)}+ Y(y)(\diff\psi+A)^2 \right]~.\eea
Note this is of the Calabi form (\ref{calabiform}). The vector field 
\bea
\frac{n+2}{n+1}\frac{\partial}{\partial\tau}\eea
is thus asymptotically the Reeb vector field: locally the metric (\ref{metric})
asymptotes, for large $\pm x$, to a cone over the local Sasaki-Einstein 
metric of \cite{paper3}.


\section{Global analysis: $\pm x > \pm x_{\pm}$}
\label{global}

We begin by making the following change of angular coordinates
\bea
\tau \, = \, -\alpha,\qquad \quad \psi \, = \, \alpha+\frac{\gamma}{n+1}~.\eea
The metric (\ref{metric}) becomes 
\bea\label{newmetric}
g  & =  & (1-x)(1-y)g_V + \frac{y-x}{4X(x)}\diff x^2 + \frac{y-x}{4Y(y)}\diff y^2 +
 \frac{v(x,y)}{(n+1)^2}\left[\diff\gamma+(n+1)A\right]^2\nonumber\\ 
&& + w(x,y)\left[\diff\alpha+\frac{f(x,y)}{n+1}\left[\diff\gamma+(n+1)A\right]\right]^2 \eea
where we have defined
\bea
w(x,y) & = & \frac{1}{y-x}\left[y^2X(x)+x^2Y(y)\right]\\
f(x,y) & = & 1-\frac{[yX(x)+xY(y)]}{y^2X(x)+x^2Y(y)}\\
v(x,y) & = & \frac{X(x)Y(y)}{w(x,y)}~.\eea
The strategy for extending the local metric (\ref{metric}) to a complete 
metric on a non-compact manifold is as follows. We shall take
\bea
y_1\leq y\leq y_2\eea
where $y_1,y_2$ are two appropriate adjacent zeroes of $Y(y)$ satisfying
\bea
y_1<y_2<1~.\eea
On the other hand, we take $x$ to be a non-compact coordinate, with 
\bea
-\infty< x\leq x_-\leq y_1~,\qquad \mathrm{or} \qquad \; 1\leq x_+ \leq x < +\infty~.
\eea
Here $x_-$ is the smallest zero of $X(x)$ and $x_+$ is the largest zero; 
thus $X(x)>0$ for all $x<x_-$ or
$x > x_+$. First, we examine 
regularity of the metric (\ref{newmetric}) for $\pm x > \pm x_{\pm}$.
Following a strategy similar to 
\cite{paper2,paper3}, we show that the induced metric  
at any constant $x$, that is not a zero of $X(x)$, 
may be extended to a complete metric on the 
total space of a $U(1)$ principal 
bundle (with local fibre coordinate $\tilde\alpha$) over
a smooth compact base space $Z(V)$. In particular, this will fix 
the metric parameter $\nu$. The analysis essentially 
carries over from that presented in \cite{paper3}. Note, however, 
that the results  of subsection \ref{lemmas} 
complete the discussion in reference \cite{paper3}. 
The remaining sections of the paper will deal with 
regularity of the metric at $x=x_{\pm}$.

\subsection{Zeroes of the metric functions}\label{zeroes}

Recall that
\bea
Y(y) \, =\,  \frac{p(y)}{(1-y)^n}
\eea
where
\bea
p(y) \, = \, (1-y)^{n+1} -\frac{(n+1)}{(n+2)}(1-y)^{n+2}-2\nu~.
\eea
One easily verifies that $p^{\prime}(y)=0$ if and only if $y=0$ or 
$y=1$. The former is a local maximum of $p(y)$, whereas the latter is 
a local minimum or a point of inflection depending on whether $n$ is 
odd or even, respectively. Defining
\bea
\numax\, =\, \frac{1}{2(n+2)}
\eea
we also see that $p(0)\leq 0$ for $\nu\geq \numax$ and $p(1)\geq 0$
for $\nu\leq 0$. Since for regularity we require two adjacent real zeroes 
$y_1,y_2$ of $Y(y)$, 
with $1\notin (y_1,y_2)$, it follows that we must take
\bea
0\leq \nu \leq \numax~.\eea
Since $p^{\prime}(0)=0$, the roots then satisfy $y_1\leq0$, $y_2\geq 0$, 
and we have $Y(y)>0$ for $y\in (y_1,y_2)$.  We also note 
that for any zero $y_i$ of $Y(y)$ we have
\bea\label{remarkable}
Y'(y_i) \, = \, -(n+1)y_i~.
\eea
The metric with $\nu=0$ is the local Ricci-flat K\"ahler metric one obtains by using
the Calabi ansatz with local K\"ahler-Einstein 
metric $g_T$ in (\ref{gtrans}). Thus this metric admits a Hamiltonian two-form 
of order one, with $S=1$. The local metric $g_T$ extends 
to a smooth K\"ahler-Einstein 
metric on a complete manifold 
only when $V=\mathbb{CP}^n$, in which case $g_T$ 
is the K\"ahler-Einstein metric on $\mathbb{CP}^{n+1}$. 
More generally, the quasi-regular Sasaki-Einstein metrics 
constructed in \cite{paper3} lead to smooth complete orbifold metrics. 
For $\nu=\numax$ one finds that $y_1=y_2=0$.
One can verify that the metric also reduces to the Calabi 
ansatz in this limit, with 
the product K\"ahler-Einstein metric on $\mathbb{CP}^1\times V$ 
as base. Hence one obtains a complete Ricci-flat K\"ahler 
metric on the canonical line bundle 
over $\mathbb{CP}^1\times V$. Henceforth we take $\nu\in (0,\nu_{\mathrm{max}})$.

We now turn to an analysis of the zeroes of $X(x)$. Recall that
\bea
X(x) \, =\,  \frac{q(x)}{(x-1)^n}
\eea
where
\bea
q(x)\,  =\,  (x-1)^{n+1} +\frac{(n+1)}{(n+2)}(x-1)^{n+2}+2\mu~.
\eea
Any root $x_0(\mu)$ of $q(x)$ therefore satisfies
\bea\label{diffx}
\frac{\diff x_0}{\diff \mu}  \, = \, - \frac{2}{(n+1)x_0(x_0-1)^n}~.
\eea 
Note that $x_0=0$ is a root of $q(x)$ when 
\bea
\mu\, =\, \bar\mu  \, =\,  \frac{(-1)^n}{2(n+2)}
\eea
and that $x_0=1$ is a root of $q(x)$ when $\mu=0$.
In our later analysis we shall require there to exist either 
a smallest zero $x_-$ of $X(x)$, with $x_-\leq y_1< 0$, or 
a largest zero $x_+$, with $1\leq x_+$. It 
is easy to see that $q'(x)=0$ if and only if $x=0$ or $x=1$. 
The former is local maximum for $n$ odd and a local minimum for 
$n$ even, while the latter is local minimum for $n$ odd 
and a point of inflection for $n$ even. The behaviour of $x_-$ is summarised 
by the following
\begin{lemma}\label{freak} 
For each $x_-\in (-\infty,0]$ there exists a unique $\mu$ such that $x_-$
is the smallest zero of $X(x)$. Moreover, $x_-(\mu)$ is monotonic.
\end{lemma}
\begin{proof}
The proof depends on the parity of $n$. For $n$ odd, $x_-\rightarrow -\infty$ 
as $\mu\rightarrow\infty$. Since $x_-(\bar{\mu})=0$, (\ref{diffx}) shows 
that $x_-(\mu)$ is monotonic decreasing in $[\bar{\mu},\infty)$. 
For $n$ even, instead $x_-\rightarrow-\infty$ as $\mu\rightarrow-\infty$. 
Equation (\ref{diffx}) now shows 
that $x_-(\mu)$ is monotonic increasing in $(-\infty,\bar{\mu}]$.
\end{proof}
For $x_+$ we similarly have
\begin{lemma}
\label{ape}
For each $x_+\in[1,\infty)$ there exists a unique $\mu\leq0$ such that $x_+$ is
the largest zero of $X(x)$. Moreover, $x_+(\mu)$ is monotonic decreasing.
\end{lemma}
\begin{proof}
Noting that $q(1)=2\mu$, and $q'(x)>0$ for $x>1$, 
we see that for a zero $x_+\geq 1$ of $X(x)$ to exist, we must
require $\mu\leq 0$ (independently of the parity of $n$). Moreover, 
(\ref{diffx}) immediately implies that $x_+(\mu)$ is monotonic 
decreasing in $\mu$.
\end{proof}

\subsection{Regularity for $\pm x > \pm x_{\pm}$}
\label{regularitysec}

Let us fix $x$ with $\pm x > \pm x_{\pm}$ and consider the positive 
definite\footnote{For $x>x_+$ the metric (\ref{newmetric}) is negative 
definite. Henceforth all metrics we write will be positive definite.} metric 
$h_x$ given by
\bea\label{base}
\mp h_x \,= \, (1-x)(1-y)g_V + \frac{y-x}{4Y(y)}\diff y^2 +  
\frac{v(x,y)}{(n+1)^2}\left[\diff\gamma+(n+1)A\right]^2 ~.\eea
Near to a root $y_i$ of $Y(y)$ we have
\bea
Y(y) \, = \, Y'(y_i)(y-y_i) + \mathcal{O}((y-y_i)^2)~.\eea
Defining
\bea
\mp R_i^2 \,= \, \frac{(y_i-x)(y-y_i)}{Y'(y_i)}\eea
for each $i=1,2$, one easily obtains, near to $R_i=0$,
\bea
h_x = \pm (x-1)(1-y_i + \mathcal{O}(R_i^2))g_V + \diff R_i^2 + R_i^2 \left[\diff\gamma+(n+1)A\right]^2+\mathcal{O}(R_i^4) ~.\eea
Fixing a point on $V$, we thus see that the metric is regular near 
to either zero provided that  we take the 
period of $\gamma$ to be $2\pi$. This ensures that $y=y_i$ is merely 
a coordinate singularity, resembling the origin of $\R^2$ in polar 
coordinates $(R_i,\gamma)$.
The one-form\label{angularform}
\bea
\diff\gamma + (n+1)A\eea
is precisely the global angular form 
on the unit circle bundle in the canonical line bundle $K_V$ over $V$. 
Indeed, recall that 
\bea
(n+1)\diff A \,= \, 2(n+1)\omega_V \, = \, \rho_V\eea
is the curvature two-form of the anti-canonical line bundle $K_V^{-1}$.
It follows that $h_x$ extends to a smooth metric on the manifold
\bea
Z(V) \, =\, 
 K_V\times_{U(1)} S^2\eea
for all $\nu$ with $0<\nu<\numax$. That is, $Z(V)$ is the total 
space of the $S^2$ bundle over $V$ obtained 
using the $U(1)$ transition functions of $K_V$, with the natural 
action of $U(1)\subset SO(3)$ on the $S^2$ fibres.

The induced metric on $\pm x > \pm x_{\pm}$ is in fact
\bea
g_x \,=\, h_x \mp w(x,y)\left[\diff\alpha+\frac{f(x,y)}{n+1}\left[\diff\gamma+(n+1)A\right]\right]^2~.\eea 
One first notes that $w(x,y)<0$ for all $x>x_+$, and $w(x,y)>0$ for all $x<x_-$, 
$y\in [y_1,y_2]$, so that $g_x$ has positive definite signature.
We then consider the one-form $\diff \alpha + B$ where
\bea
B \,=\, \frac{f(x,y)}{n+1}\left[\diff\gamma + (n+1)A\right]~.\eea
As in \cite{paper3}, the strategy now is to show that, for appropriate 
$\nu$, there exists 
$\ell\in\R_+$ such that $\ell^{-1}B$ is locally a connection one-form on a 
$U(1)$ principal bundle over $Z(V)$. By periodically identifying $\alpha$ 
with period $2\pi\ell$, we thereby obtain a complete metric on the total space 
of this $U(1)$ principal bundle. 

Assuming that $V$ is simply-connected\footnote{Again, $b_1(V)=0$ necessarily.}, one can show that 
$Z(V)$ has no torsion in $H^2(Z(V);\Z)$. The isomorphism class of a complex 
line bundle over $Z(V)$ is then determined completely 
by the integral of a curvature 
two-form over a basis of two-cycles. Such a basis is 
provided by $\{\Sigma,\sigma_*(\Sigma_i)\}$. Here 
$\Sigma$ is represented by 
a copy of the fibre $S^2$ at any fixed point on $V$; 
$\sigma:V\rightarrow Z(V)$ is the section $y=y_1$, and $\{\Sigma_i\}$ 
are a basis of two-cycles for $H^2(V;\Z)$, which similarly is 
torsion-free. 

One may then compute the periods
\bea
\int_{\Sigma}\frac{\diff B}{2\pi} & = & \frac{1}{n+1}(f(x,y_1)-f(x,y_2))\nonumber\\
\int_{\sigma_*(\Sigma_i)}\frac{\diff B}{2\pi} & = & \frac{f(x,y_1)}{n+1}
\langle c_1(K_V^{-1}),\Sigma_i\rangle~.\eea
We now note that
\bea
f(x,y_i) \, =\,  \frac{y_i-1}{y_i}\,  \equiv \, f(y_i)
\eea
is independent of the choice of $x$. Notice 
that $f(y_1)>0$ since $y_1<0$. Defining
\bea
\ell \, = \, \frac{I f(y_1)}{k(n+1)}\eea
it follows that the periods of 
$\ell^{-1}\diff B/2\pi$ over the two-cycles $\{\Sigma,\sigma_*(\Sigma_i)\}$ are 
\bea\label{periodsone}
\left\{\frac{k(f(y_1)-f(y_2))}{If(y_1)}, \frac{k}{I}\langle c_1(K_V^{-1}),\Sigma_i\rangle\right\}~.\eea
We now choose $\nu$ so that
\bea\label{rational}
\frac{f(y_1)-f(y_2)}{f(y_1)} \, = \, \frac{pI}{k}\in\mathbb{Q}\eea
is rational. In particular $p,k\in \mathbb{N}$ are positive 
integers; no coprime condition is assumed, so the rational 
number (\ref{rational}) is not assumed to be expressed in lowest terms. 
We note the following useful identities 
\bea
\frac{y_2(1-y_1)}{y_2-y_1} & = & \frac{k}{pI}\nn\\
\frac{y_1(1-y_2)}{y_2-y_1} & = & \frac{k}{pI}-1~,
\label{bundletwists}
\eea
which we will  use repeatedly in the remainder of the paper.
We shall return to which values of $p$ and $k$ are allowed momentarily.
The periods (\ref{periodsone}) are then
\bea\label{periods}
\left\{p,k\left\langle c_1(K_V^{-1/I}),\Sigma_i\right\rangle\right\}~.\eea
Notice that, by definition of the Fano index $I$, 
$c_1(K_V^{-1/I})\in H^2(V;\Z)$ is primitive. 
Defining the new angular variable 
\bea
\tilde \alpha \, = \, \ell^{-1} \alpha ~,
\eea
it follows that if we periodically identify $\tilde\alpha$ with period $2\pi$,
then the metric $g_x$ at fixed $x$ is a smooth 
complete metric on a $U(1)$ principal bundle over $Z(V)$, where 
the Chern numbers of the circle fibration over the two-cycles 
$\{\Sigma,\sigma_*(\Sigma_i)\}$ are given by (\ref{periods}). 
We denote the total space by $L^{p,k}(V)$. Note that $L^{ph,kh}(V)$ is simply a $\Z_h$ quotiet 
of $L^{p,k}(V)$, where $\Z_h\subset U(1)$ acts on the fibres 
of the $U(1)$ principal bundle $L^{p,k}(V)\rightarrow Z(V)$.
We may also think of this manifold as a Lens space 
$L(1,p)=S^3/\Z_p$ fibration over $V$.
Since the regularity analysis was 
essentially independent of $x$, we see that the Ricci-flat K\"ahler metric $g$ 
extends to a smooth asymptotically conical 
metric on $\R_+\times L^{p,k}(V)$, where $x-x_+>0$ or $x_- -x>0$ is a coordinate on 
$\R_+$. The asymptotic cone has Sasaki-Einstein base $Y^{p,k}(V)$, 
constructed originally in \cite{paper3}.

\subsection{Allowed values of $p$ and $k$}
\label{lemmas}

We will now determine the allowed values of $p$ and $k$ in 
(\ref{rational}). We begin by 
defining
\bea
Q(\nu) \, \equiv \, \frac{f(y_1)-f(y_2)}{f(y_1)} \, = \, \frac{y_2-y_1}{y_2(1-y_1)}~,\eea
regarding the roots $y_i$ of $p(y)$ as functions of the metric parameter $\nu$.
In the remainder of this section we shall prove
\begin{proposition}\label{littleprop} The function $Q:[0,\numax]\rightarrow\R$ is 
a continuous monotonic increasing function with $Q(0)=1$, 
$Q(\numax)=2$. \end{proposition}
Given (\ref{rational}), this implies that 
\bea
\frac{pI}{2}~<~k~<~pI
\label{pkrange}
\eea
and that, for each $p$ and $k$, there is a corresponding \emph{unique} 
metric\footnote{Note this analysis completes
the argument presented in \cite{paper3}.}. From the defining equations of the roots we have 
\bea
2\nu ~=~ (1-y_1)^{n+1} - \frac{n+1}{n+2} (1-y_1)^{n+2}~ = ~(1-y_2)^{n+1} - \frac{n+1}{n+2} (1-y_2)^{n+2}~,
\label{rooots}
\eea
and from (\ref{rooots}) one easily obtains  the following useful identity
\bea
\frac{(1-y_2)^{n+1} }{(1-y_1)^{n+1} } & = &\frac{1+(n+1)y_1}{1+(n+1)y_2} ~.
\label{banana}
\eea
From (\ref{rooots}) one also computes 
\bea
\frac{\diff y_i }{\diff \nu } & = & - \frac{2}{(n+1) y_i (1-y_i )^n}~.
\label{rogerabbit}
\eea
One may then use this formula to prove the following
\begin{lemma}\label{carrot} $y_1(\nu)$ (respectively $y_2(\nu)$) is monotonically increasing (respectively decreasing)
in the interval $[0,\nu_{\mathrm{max}}]$. In particular, in the open interval $(0,\nu_{\mathrm{max}})$ the following bounds hold:
\bea
-\frac{1}{n+1} <& y_1& < 0\\
0 < & y_2 & < 1~.
\eea
\end{lemma}
\begin{proof} We have $y_1(0)=-1/(n+1)$ and $y_1(\nu)<0$ 
for all $\nu\in (0,\numax)$ since $p^{\prime}(y=0)=0$ for all $\nu$. Thus from (\ref{rogerabbit}) $y_1$ 
is monotonic increasing in this range. 
A similar argument applies for $y_2$ on noting that $y_2>0$ and $y_2(0)=1$.
\end{proof}
We now define 
\bea
R(\nu) &=& \frac{y_1(1-y_2)}{y_2(1-y_1)}~,
\eea
so that $Q=1-R$. Using (\ref{rogerabbit}) and (\ref{banana}) one easily obtains
\bea
\frac{\diff R}{\diff \nu}  ~= ~\frac{-2(y_2-y_1)}{(n+1)y_1y_2^3(1-y_1)(1-y_2)^n(1+(n+1)y_2)}D
\eea
where we have defined
\bea
D \, \equiv \,  y_1+ y_2 + (n+1)y_1 y_2~.
\eea
Making use of the above identities, we also compute
\bea
\frac{\diff D}{\diff \nu} \, = \, - \frac{2(1+(n+1)y_2)}{(n+1)y_1(1-y_1)^n} (1+R)~.
\eea
The above computations, and 
Lemma \ref{carrot}, result in
\begin{lemma}\label{hedgehog} In the open interval $(0,\nu_{\mathrm{max}})$ 
the sign of $D'$ is correlated with that of $1+R$, and the 
sign of $R'$ is correlated with that of $D$. In particular, we have
\bea
R ~=~ -1 & \mathrm{iff} &  \frac{\diff D}{\diff \nu} = 0 \\
D~ =~ 0 & \mathrm{iff}  & \frac{\diff R}{\diff \nu} = 0 ~.
\eea\end{lemma}

Next, we turn to analysing the behaviour of $R$, $D$ and their derivatives 
at the endpoints 
of the interval.  It is easily checked that $R(0)=0$, $D(0)=-1/(n+1)$ and $D(\nu_{\mathrm{max}})=0$.
In order to compute $R(\numax)$ we write $\nu=\nu_{\mathrm{max}}- (n+1)\delta^2$, 
where the factor of $(n+1)$ is inserted for later convenience. 
We may solve for $y$ in a power series in $\delta$; the first two terms
suffice for our purposes:
\bea
y_1 ~= ~-2\delta +\frac{4}{3}n\delta^2  +{\cal O}(\delta^{3})~,\qquad y_2 ~= ~ 2\delta 
+\frac{4}{3}n\delta^2 + {\cal O}(\delta^{3})~.
\label{rambo}
\eea
With these one then computes
\bea
R \, =\,  -1 +\left(4+\frac{4}{3}n\right)\delta + {\cal O}(\delta^2)~,
\label{pile}
\eea
which proves that $R\to -1$ as $\nu\to\nu_{\mathrm{max}}$.
Before turning to the proof of Proposition \ref{littleprop}
we shall
need another 
\begin{lemma}\label{hamster} $ R'(\nu) \to -\infty$ for $\nu\to 0^+$ and $\nu\to \nu_{\mathrm{max}}^-$. \end{lemma}
\begin{proof} Near $\nu=0$ this is easily checked; near $\nu=\nu_{\mathrm{max}}$ the result follows from (\ref{pile}).\end{proof}

\begin{proof} of Proposition \ref{littleprop}. It is enough to show that
$R(\nu)$ is a monotonically decreasing function in the interval $(0,\nu_{\mathrm{max}})$. Suppose this is not so. Then $R'(\nu_1)=0$ for some least 
$\nu_1\in (0,\numax)$. By Lemma \ref{hedgehog}, this is also the first 
time that $D$ crosses zero, $D(\nu_1)=0$. There are then three cases. 
We make repeated use of Lemma \ref{hedgehog}:
\begin{itemize}
\item Suppose $R(\nu_1)>-1$. Then $D'(\nu_1) > 0$, and $D$ is positive 
in the range $(\nu_1,\nu_1+\epsilon)$ for some $\epsilon>0$. 
Since $D(\numax) = 0$, there must be a turning point $D'(\nu_2) = 0$ 
for some smallest $\nu_2\in (\nu_1,\numax)$. We then 
have $R(\nu_2) = -1$. But $R'(\nu) > 0$ for all $\nu\in (\nu_1,\nu_2)$ 
since $D>0$ in this interval, which implies that $R$
is strictly monotonic increasing in this range. 
This is a contradiction since we assumed $R(\nu_1) > -1$.
\item Suppose $R(\nu_1) < -1$. Then $D'(\nu_1) < 0$. This is an immediate
contradiction, since $D(0)=-1/(n+1)$ is negative and $\nu_1$ 
is the first zero of $D$; hence $D'(\nu_1)$ must be non-negative.
\item Suppose $R(\nu_1)=-1$. Either 
$R(\nu_1)=-1$ is a local minimum or a point of inflection:
\begin{itemize}
\item Suppose $R(\nu_1)=-1$ is a local minimum of $R$. Since $R(\numax)=-1$ 
also, there must be a turning point $R'(\nu_3)=0$ for some least 
$\nu_3\in (\nu_1,\numax)$. Then 
$D(\nu_1)=D(\nu_3)=0$. Since $R(\nu)>-1$ for $\nu\in (\nu_1,\nu_3)$, it 
follows that $D'>0$ in the same range, a contradiction.
\item Suppose $R(\nu_1)=-1$ is a point of inflection of $R$. Then 
$D(\nu_1)=0$ is a local maximum of $D$. But since $D(\numax)=0$ also, 
$D$ must have a turning point $D'(\nu_5)=0$ for some least $\nu_5\in (\nu_1, 
\numax)$. Thus $R(\nu_5)=-1$. But $D(\nu)<0$ 
for all $\nu\in (\nu_1,\nu_5)$ implies that $R'<0$ in the same range, a
contradiction.
\end{itemize}
\end{itemize}
This proves that there is no $\nu_1\in (0,\numax)$ where $R'(\nu_1)=0$, and 
hence $R$ is strictly monotonic decreasing in this range.
\end{proof}

\subsection{Summary}
 
We end the section by summarising what we have proven so far:
\begin{proposition} Let $(V,g_V)$ be a complete 
K\"ahler-Einstein manifold of positive 
Ricci curvature with Fano index $I$. 
Then for every $p, k\in \mathbb{N}$
positive integers with $pI/2 < k < pI$ there is an 
asymptotically conical Ricci-flat K\"ahler metric on $\R_+\times L^{p,k}(V)$, 
with 
local form (\ref{metric}) - (\ref{functions}) and $\R_+$ coordinate 
either $x-x_+>0$ or $x_--x>0$. 
All metric parameters and ranges of coordinates are fixed uniquely
for a given $p$ and $k$, except for the constant $\mu$. The metric is 
asymptotically a cone over the Sasaki-Einstein manifold $Y^{p,k}(V)$.
\end{proposition}

In the remainder of the paper we examine regularity of the above metric 
at $x=x_{\pm}$ for the cases $x_-=y_1$, $x_+=1$ and
$x_-< y_1$, $x_+> 1$.


\section{Small resolutions}
\label{smallsection}

In this section we consider the cases $x_-=y_1$ and $x_+=1$. 
Equivalently, these are the special 
cases $\mu=\pm\nu$ and $\mu=0$. These will give
rise to  \emph{partial} small resolutions, where one blows up the 
Fano $V$ or a weighted projective space
$\mathbb{WCP}^1$, respectively. In particular, we prove 
Theorem \ref{orbifolds} and Theorem \ref{flopper}.
As a simple consequence of Theorem \ref{orbifolds}, we shall also 
obtain \emph{smooth} small resolutions, as
summarised in Corollary \ref{p=1}. The remaining cases 
where $x_-<y_1$ or $x_+>1$
will be the subject of section \ref{bigones}.

\subsection{Partial resolutions I: $x_-=y_1$}

In this section we analyse regularity in the case that $x_-=y_1$. 
This special case arises since the function $y-x$ appears in 
the metric (\ref{metric}); when $x_-<y_1$ this function is strictly 
everywhere positive, whereas when $x_-=y_1$ the function has 
a vanishing locus. From section \ref{zeroes}
one easily deduces that $x_-=y_1$ corresponds to $\mu=-\nu$ when 
$n$ is odd, and $\mu=\nu$ when $n$ is even.

The first remarkable point to note 
is that, due to the symmetry in 
$x$ and $y$, the analysis of the collapse at $x=x_-$, for $y_1<y<y_2$, 
is identical to that of the collapse at $y=y_i$ for $x<x_-$. 
Thus for fixed $y\in (y_1,y_2)$ we deduce immediately that the metric 
collapses smoothly at $x=x_-$ for all $\mu$. It thus remains to check 
the behaviour of the metric at $\{x=x_-,y=y_1\}$ and 
$\{x=x_-,y=y_2\}$. 

The induced metric on $x=y_1$ is 
\bea
g_{\{x=y_1\}} =  (1-y_1)(1-y)g_V \!+\! \frac{\diff y^2}{4W(y)} 
+ y_1^2 W(y)\left[\diff\alpha\!+\!\frac{y_1-1}{(n+1)y_1}[\diff\gamma+(n+1)A]\right]^2 .
\label{bob}
\eea
where we have defined
\bea
W(y) \, = \, \frac{Y(y)}{y-y_1}~.
\eea
Since $y=y_1$ is a simple zero of $Y(y)$, for $y_1\neq0$, we have 
$W(y_1)\neq 0$. Thus $y=y_2$ is the only zero of $W(y)$ for 
$y_1\leq y\leq y_2$, with $W(y)>0$ for $y_1\leq y<y_2$. We thus 
see that the above metric is regular at $y=y_1$. 

In general it will turn out that $\{x=y_1, y=y_2\}$ is a locus 
of orbifold singularities. Although one can analyse the behaviour of 
the metric here using (\ref{bob}), it will turn out that $x=y_1$ 
is a rather unusual type of coordinate 
singularity; this might have been anticipated from the above 
factorisation of $Y(y)$ into $W(y)$. We will therefore introduce 
a new set of coordinates, resembling polar coordinates\footnote{As we 
shall see later, the global structure is in fact $\C^2/\Z_p$.} on 
$\R^4=\R^2\oplus\R^2$, in which this coordinate singularity is more 
easily understood.

We begin by defining 
\bea
R_1^2 & = & a_1 (x-y_2)(y-y_2)\nn\\
R_2^2 & = & a_2 (x-y_1)(y-y_1)
\label{newpolar}
\eea
where $a_1,a_2$ are constants that will be fixed later. 
The induced change in the $x-y$ part of the metric is then
\bea
&&\!\!\!\!\!\frac{y-x}{4X(x)}\diff x^2 + \frac{y-x}{4Y(y)}\diff y^2 = \\
&& \!\!\!\!\!\!\!\!\!\!\!\!\! = \frac{1}{(y-x)(y_1-y_2)^2} 
\Bigg\{\! \left[ \frac{(x-y_1)^2}{X(x)}\! +\! \frac{(y-y_1)^2}{Y(y)}\right] \frac{R_1^2}{a_1^2}\diff R_1^2 
+\!\! \left[ \frac{(x-y_2)^2}{X(x)} \! +\! \frac{(y-y_2)^2}{Y(y)}\right] \frac{R_2^2}{a_2^2}\diff R_2^2\nn\\
&& \;\;\;\;\;\qquad\qquad\qquad- \frac{2R_1 R_2}{a_1a_2} \Bigg[ \frac{(x-y_1)(x-y_2)}{X(x)} 
+ \frac{(y-y_1)(y-y_2)}{Y(y)}\Bigg]
\diff R_1 \diff R_2 \Bigg\}~.\nn \eea
Let us expand this near to $x=y_1$, $y=y_2$, which is $R_2=0$, $R_1=0$. Using 
\bea
y_2-y & = & \frac{R_1^2}{a_1(y_2-y_1)}+\mathcal{O}(R_1^2, R_2^2)\nn\\
y_1 - x & = & -\frac{R_2^2}{a_2(y_2-y_1)} + \mathcal{O}(R_1^2, R_2^2)
\eea
one finds
\bea
&&\quad\frac{y-x}{4X(x)}\diff x^2 + \frac{y-x}{4Y(y)}\diff y^2 =\\
&&= -\frac{1}{a_1 Y'(y_2)} [1+\mathcal{O}(R^2)]\diff R_1^2  
+ \frac{1}{a_2 X'(y_1)} [1+\mathcal{O}(R^2)]\diff R_2^2 + 
\mathcal{O}(R_1,R_2)\diff R_1\diff R_2~,\nn
\eea
where $\mathcal{O}(R^2)$ denotes terms of order $\mathcal{O}(R_1^2)$ or order $\mathcal{O}(R_2^2)$. We thus set
\bea
a_1  =  - \frac{1}{Y'(y_2)} = \frac{1}{(n+1)y_2} \qquad a_2  =  \frac{1}{X'(y_1)} = \frac{1}{(n+1)y_1}~.
\eea
The change of coordinates (\ref{newpolar}) now becomes
\bea\label{newnewpolar}
R_1^2 &=& \frac{1}{(n+1)y_2}(y_2-x)(y_2-y)\nn\\
R_2^2 &= &-\frac{1}{(n+1)y_1}(y_1-x)(y-y_1)\eea
and these relations imply
\bea
\frac{y_2}{y_2-x} R_1^2 +\frac{y_1}{x-y_1} R_2^2 & = & \frac{y_2-y_1}{n+1} \nn\\
\frac{y_2}{y_2-y} R_1^2 + \frac{y_1}{y-y_1} R_2^2 & = & \frac{y_2-y_1}{n+1}~. 
\eea
\begin{figure}[!ht]
\vspace{5mm}
\begin{center}
\epsfig{file=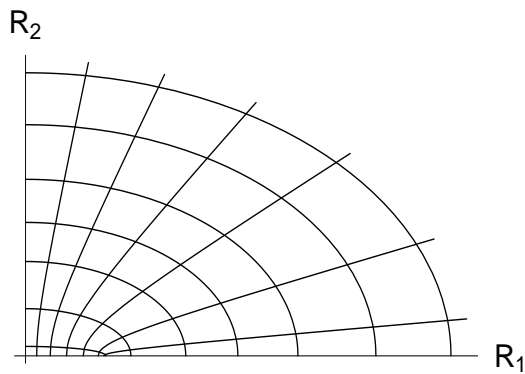,width=7.5cm,height=5cm}
\end{center}
\caption{An illustration of the change of coordinates (\ref{newnewpolar}). Curves of constant
$x$ and constant $y$ are depicted.} 
\vspace{5mm}
\label{pictu}
\end{figure}
Despite the symmetry in $x$ and $y$, the curves of constant $x$ and constant $y$ are different due 
to the difference in ranges of the variables. Recall that $y_1\leq y\leq y_2$ and 
$x\leq x_-=y_1$. The constant $x$ curves are ellipses, while the constant $y$ curves are
hyperbolae. This behaviour is depicted in Figure 1.
Notice that both sets of curves degenerate on the $R_1$-axis.
Indeed, note that $R_1=0$ if and only if $y=y_2$; but $R_2=0$ if $x=y_1$ \emph{or}
$y=y_1$. In particular, when $x=y_1$ we have 
$R_1^2=a_1(y_2-y)(y_2-y_1)$ and thus this branch of the $R_1$-axis is coordinatised 
by $y$. On the other hand, when $y=y_1$ we have $R_1^2 = a_1(y_2-x)(y_2-y_1)$, 
and thus this branch of the $R_1$-axis is coordinatised by $x$. 

Fixing a point on $V$, to leading order the induced metric near $R_1=R_2=0$ is
\bea\label{fibreguy}
g_{\mathrm{fibre}} & = & \diff R_1^2 + \frac{(n+1)^2}{(y_2-y_1)^2}
y_2^2 R_1^2 \left[y_1 \diff \alpha + \frac{(y_1-1)}{(n+1)} \diff \gamma \right]^2  \nn\\
&+& \diff R_2^2 +\frac{(n+1)^2}{(y_2-y_1)^2} y_1^2 R_2^2 \left[y_2 \diff \alpha + \frac{(y_2-1)}{(n+1)} \diff \gamma \right]^2~.
\eea
We then define
\bea
\phi_1 & = & -\frac{(n+1)}{(y_2-y_1)} y_2 \left[y_1  \alpha + \frac{(y_1-1)}{(n+1)} \gamma\right]\nn\\
\phi_2 & = & -\frac{(n+1)}{(y_2-y_1)} y_1  \left[y_2  \alpha + \frac{(y_2-1)}{(n+1)} \gamma\right]
\label{polarbears}
\eea
so that (\ref{fibreguy}) becomes
\bea\label{flatboy}
g_{\mathrm{fibre}} = \diff R_1^2 + R_1^2 \diff\phi_1^2 + \diff R_2^2 + R_2^2 \diff\phi_1^2~.\eea
In terms of the variable $\tilde\alpha =\ell^{-1}\alpha$, the change of coordinates 
(\ref{polarbears}) becomes
\bea\label{transformers}
\phi_1 & = & \frac{1}{p}\tilde{\alpha} + \frac{k}{pI}\gamma\nn\\
\phi_2 & = & \frac{1}{p}\tilde{\alpha} + \left(\frac{k}{pI}-1\right)\gamma\eea
on using the identities (\ref{bundletwists}). 
Notice that the Jacobian of the transformation (\ref{transformers}) is $1/p$. 
Recall also from section \ref{regularitysec} that $\tilde\alpha$ and 
$\gamma$ are periodically identified with
period $2\pi$. It follows from (\ref{flatboy}) that a neighbourhood 
of $R_1=R_2=0$, at a fixed point on $V$, is diffeomorphic to 
$\R^4/\Z_p$. Indeed, as mentioned in section \ref{global}, 
the surfaces of constant $x<x_-$ are Lens space fibrations $L(1,p)=S^3/\Z_p$ 
over $V$. These are then constant radius surfaces in the 
$\R^4/\Z_p$ fibration over $V$. The set of points
$\{x=x_-,y=y_2\}$ are the zero-section, which is a copy of $V$ and 
locus of orbifold singularities. In fact, the 
possible existence of such metrics was raised at the end of 
reference \cite{paper2}. The fibres must of course be complex 
submanifolds, and one easily checks that the complex structure is 
such that each fibre is $\C^2/\Z_p$.

One may work out the precise fibration structure as follows. 
Setting $y=y_1$ and $y=y_2$ gives two different $\C/\Z_p$ fibrations 
over $V$. It is enough to determine the fibration structure of these 
bundles. Of course, the unit circle bundle in each is a 
$U(1)$ principal bundle over $V$. These $U(1)$ bundles are determined from 
the analysis in section \ref{global}. The $U(1)$ bundle 
at fixed $x<x_-$ and $y=y_1$ has first Chern class 
$k c_1(K_V)/I$. The associated complex line bundle is 
thus $K_V^{k/I}$. The $U(1)$ bundle at fixed $x<x_-$ and 
$y=y_2$ is determined from the periods of 
$-\ell^{-1}\diff B/2\pi$ over the image of 
cycles $\Sigma$ in $V$ at $y=y_2$. We denote these as 
$\tau_*(\Sigma_i)$. The periods are given by 
\bea
-\int_{\tau_*(\Sigma_i)}\ell^{-1}\frac{\diff B}{2\pi} & = 
& -\ell^{-1}\frac{f(y_2)}{n+1}
\langle c_1(K_V^{-1}),\Sigma_i\rangle\nn \\
& = & \left(p-\frac{k}{I}\right)\langle c_1(K_V^{-1}), \Sigma_i\rangle~.\eea
This implies that the $U(1)$ bundle has first 
Chern class $(pI-k)c_1(K_V)/I$, and thus the associated line 
bundle is $K_V^{(pI-k)/I}$.

It is now a simple matter to determine the twisting of the 
$\C^2/\Z_p$ fibres themselves. The $\Z_p\subset U(1)\subset SU(2)$ 
acts via the standard action of $SU(2)$ on $\C^2$. Define $\mathcal{L}_1=K_V$ 
and $\mathcal{L}_2=K_V^{k/I}$. Let $\{U_{\alpha}\}$ be a trivialising 
open cover of $V$, and let $g^i_{\alpha\beta}$, $i=1,2$, denote the 
transition functions of the above bundles. Thus $g^i_{\alpha\beta}:
U_{\alpha}\cap U_{\beta}\rightarrow S^1$. 
Let $(z_1,z_2)$ denote standard complex coordinates on $\C^2$. 
These are identified via the action of $\Z_p\subset U(1)\subset SU(2)$. 
We must specify precisely how $U_{\alpha}\times \C^2/\Z_p$ is 
glued to $U_{\beta}\times \C^2/\Z_p$ over the overlap 
$U_{\alpha\beta}= U_{\alpha}\cap U_{\beta}$. To do this, we define 
the following action of $\T^2=S^1\times S^1$ on $\C^2/\Z_p$:
\bea\label{actionman}
\lambda &: & S^1 \times S^1 \times \C^2/\Z_p  \rightarrow \C^2/\Z_p\nonumber\\
&& \lambda(\theta_1,\theta_2;z_1,z_2) = (\exp(i\theta_1-i\theta_2/p)z_1, 
\exp(i\theta_2/p)z_2)~.
\eea
Note that this indeed defines an action of $S^1\times S^1$ on 
$\C^2/\Z_p$. Note also that the standard action of $U(1)\subset SU(2)$ on 
$\C^2$ descends to a non-effective action of $U(1)$ on the quotient
$\C^2/\Z_p$ -- this factors $p$ times through the effective $U(1)$ 
action in (\ref{actionman}). 
The $\C^2/\Z_p$ bundle is then constructed using the gluing functions
\bea
F_{\alpha\beta} & : & U_{\alpha\beta}\times\C^2/\Z_p\rightarrow U_{\alpha\beta}\times\C^2/\Z_p\nonumber\\
&& F_{\alpha\beta}[u;z_1,z_2] = [u;\lambda(g_{\alpha\beta}^1(u),g_{\alpha\beta}^2(u); z_1,z_2)]~.\eea

To check this is correct, we simply set $z_1=0$ and $z_2=0$ separately. 
This should be equivalent to setting $y=y_1$ and $y=y_2$, respectively, 
to give $\C/\Z_p$ fibrations over $V$. From (\ref{actionman}) we 
see that $z_1=0$ has $U(1)$ principal bundle given by 
$K_V^{k/I}$. On the other hand, setting $z_2=0$, the corresponding $U(1)$ 
principal bundle is given by $K_V^p\otimes K_V^{-k/I}$. These 
are precisely the same $\C/\Z_p$ fibrations determined above using the 
metric. This completes the proof of Theorem \ref{orbifolds}.

\subsection{Smooth resolutions: $p=1$}

Setting $p=1$ in the last subsection gives a family of 
smooth complete Ricci-flat K\"ahler metrics for each choice of $(V,g_V)$. 
These are all holomorphic $\C^2$ fibrations over $V$. From (\ref{actionman}) 
this is easily seen to be a direct sum of two complex line bundles 
over $V$, namely $\left[K_V\otimes K_V^{-k/I}\right]\oplus K_V^{k/I}$. 
Setting $m=I-k$ this is $K_V^{m/I}\oplus K_V^{(I-m)/I}$, as 
stated in Corollary \ref{p=1}. The range of $k$ is given by 
(\ref{pkrange}) with $p=1$, which implies that $0<m<I/2$.

For example, we may take $V=\mathbb{CP}^n$ with its standard K\"ahler-Einstein
metric. In this case $I=n+1$ and $K_V=\mathcal{O}(-(n+1))$, so that 
$K_V^{-1/I}=\mathcal{O}(1)$. We have $0<m<(n+1)/2$ and the 
metrics are defined on the total space of the rank two holomorphic 
vector bundle  $\mathcal{O}(-m)\oplus\mathcal{O}(-(n-m+1))$ over 
$\mathbb{CP}^n$. Note that $m=n=1$ is the 
small resolution of the conifold ${\cal O}(-1)\oplus {\cal
O}(-1)\to \mathbb{CP}^1$, which is understood as a limiting case. 

As another simple example, one might take $V$ to be a product of 
complex projective spaces, equipped with the natural product K\"ahler-Einstein
metric:
\bea\label{product}
V = \prod_{a=1}^M \mathbb{CP}^{d_a-1}\eea
where
\bea
\sum_{a=1}^M d_a = n+M~.\eea
In this case $I=\mathrm{hcf}\{d_a\}$ and the rank two holomorphic 
vector bundle is given by
\bea
\mathcal{O}(-md_1/I,\ldots,-md_M/I)\oplus\mathcal{O}(-(I-m)d_1/I,\ldots,-(I-m)d_M/I)~.\eea

\subsection{Partial resolutions II: $x_+=1$}
\label{floppy}

In this subsection we consider 
a different small partial resolution, where the Fano $V$
shrinks to zero size, 
while a weighted projective space $\mathbb{WCP}^1$ is blown up. 
The analysis of this subsection is summarised 
by Theorem \ref{flopper}. Notice that when $V=\mathbb{CP}^1$, 
the process of going from 
the round $\mathbb{CP}^1$ to the weighted $\mathbb{WCP}^1$ may be 
understood as a form of generalised flop transition (called a \emph{flip})
in the K\"ahler moduli space of the family of toric Gorenstein 
singularities $C(Y^{p,q})$  \cite{toricpaper}. This will be discussed elsewhere \cite{np2}.

From the general form of the metric (\ref{metric}) it is simple to see 
that in order for $V$ to collapse
one requires $x_+=1$, and this implies that $\mu=0$. We then have 
\bea
X(x) \, = \, \frac{x-1}{n+2}(1+(n+1)x)~.
\eea 
We choose the $x\geq 1$ branch of $x$.
To analyse the metric near $x=1$ it is useful to 
change coordinates, defining
\bea
x \, = \, 1+ r^2~. 
\eea
Expanding the metric near $r=0$, and keeping terms up to order $r^2$, we find 
\bea\label{longmetric}
g & = & (1-y) \Bigg\{ \diff r^2 + r^2\Big[ g_V + \frac{1}{(n+1)^2} \left[ \diff \gamma +
 (n+1)(A + F(y)\diff\alpha)\right]^2\nn\\
&+& \frac{\diff y^2}{4Y(y)(1-y)} + G(y)\diff\alpha^2 \Big]\Bigg\}
 + \frac{1-y}{4Y(y)}\diff y^2 + \frac{Y(y)}{1-y}\diff \alpha^2 +\mathcal{O}(r^4)~,
\eea
where we have defined
\bea
F(y) &=& \frac{Y(y)}{(1-y)^2}-\frac{y}{1-y} \\
G(y) &=& \frac{1}{(1-y)^4}\left[Y(y)^2+(1-y)-2Y(y)(1-y)\right]~.\eea

We first analyse the induced metric on $r=0$, which is given by
\bea
g_W & = &\frac{1-y}{4Y(y)}\diff y^2 + \frac{Y(y)}{1-y}\diff \alpha^2~.\eea
As usual, near each root $y_i$ we introduce the coordinates
\bea
R^2_i \,= \, \frac{(y_i-1)}{(n+1)y_i}(y-y_i)
\eea
from which we see that
\bea
g_W &=&  \diff R^2_i + \left[\frac{I}{k}\frac{y_i(1-y_1)R_i}{y_1(1-y_i)}\right]^2
\diff \tilde \alpha^2+\mathcal{O}(R_i^4)~.
\eea
Thus 
\bea
g_W & = & \left\{
\begin{array}{l}
  \diff R^2_1 + \frac{I^2}{k^2}R^2_1\diff \tilde\alpha^2 +\mathcal{O}(R_1^4)\quad \mathrm{near} \quad y= y_1\nn\\[2mm]
  \diff R^2_2 + \frac{I^2}{(pI-k)^2}R^2_2\diff \tilde\alpha^2 +\mathcal{O}(R_2^4)\quad \mathrm{near} \quad y= y_2~.
\end{array}\right.
\eea
Recall now that $\tilde\alpha$ has period $2\pi$. In order to obtain 
an orbifold singularity near to $y=y_1$, it is therefore necessary that 
the Fano index $I$ divides $k$. Thus we assume this, and define 
$k=Id$ with $d$ a positive integer. It follows that $g_W$ smoothly 
approaches, in an orbifold sense, the flat metric on $\C/\Z_d$, 
where $y=y_1$ is the origin. Similarly, at $y=y_2$ the metric smoothly 
approaches the flat metric on $\C/\Z_{p-d}$. It follows that, 
provided $k=Id$, the induced metric on $r=0$ is a smooth K\"ahler 
orbifold metric on the weighted projective space $W=\mathbb{WCP}^1_{[d,p-d]}$.

Now fix any smooth point $(y,\alpha)$ on $W$, so $y_1<y<y_2$. Setting 
$R^2 = (1-y)r^2$, the induced metric near to $R=0$ is
\bea\label{rudedude}
g &=& \diff R^2 + R^2\left[ g_V + \frac{1}{(n+1)^2}\left[
\diff\gamma+(n+1)A\right]^2\right]+\mathcal{O}(R^4)~.\eea
From section \ref{global}, $\gamma$ has period $2\pi$, and 
the induced metric (\ref{rudedude}) is simply the canonical complex 
cone $\C_V$ over $V$. Equivalently, fixed $R$ gives the associated circle bundle 
in the 
canonical line bundle over $(V,g_V)$, and near to $R=0$ 
the whole metric is a real cone 
over this regular Sasaki-Einstein manifold. Thus near to $R=0$ the 
fibre metric over a smooth point on $W$ 
itself approaches a Ricci-flat K\"ahler cone.

One needs to consider what happens over 
the roots $y=y_i$ separately. These are the 
singular points of the weighted projective 
space $W=\mathbb{WCP}^1_{[d,p-d]}$. 
To determine the period of $\gamma$ in (\ref{rudedude}) 
over these points one may simply compute the volume of 
$\{y=y_i,x=1+r^2\}$ with $r^2$ small and fixed, 
and compare with section \ref{global}. 
Each space is a $U(1)$ principal bundle over $V$, namely that 
associated to the complex line bundles $K_V^{k/I}$, $K_V^{(pI-k)/I}$, 
respectively. From section \ref{global} we have the induced metric
\bea
g \, = \, (1-y_i)r^2g_V + \frac{y_i^2}{1-y_i}r^2\diff\alpha^2 +\mathcal{O}(r^4)
\eea
where $\alpha$ has period $2\pi\ell$. Comparing with (\ref{rudedude}), 
we see that $\gamma$ must be identified with period $\Delta\gamma_i$ 
over each pole, where
\bea
\Delta\gamma_i \, = \, \frac{(n+1)y_i\ell}{1-y_i}2\pi \, = \, 
 \left\{
 \begin{array}{l}
  -\frac{2\pi I}{k} \quad i=1\\[2mm]
  \frac{2\pi I}{pI-k} \quad i=2~.
\end{array}\right.\eea
Thus $\gamma$ has period $2\pi/d$ over $y=y_1$ and period $2\pi/(p-d)$ over 
$y=y_2$. This implies that the fibres over the singular points of 
$W$ are complex cones 
over $V$ associated to $K_V^{d}$ and $K_V^{p-d}$, respectively; 
the generic fibre is the complex cone $\C_V$ associated to $K_V$. 
This gives fibres $\C_V/\Z_d$ and $\C_V/\Z_{p-d}$ over the singular 
points of $W$, respectively.

In fact this latter behaviour of the fibres 
could have been deduced differently, by considering the 
fibration structure. Recall we have now checked that the 
metric is smooth away from $W=\mathbb{WCP}^1_{[d,p-d]}$, 
that $W$ is itself a smooth orbifold, and that each fibre over $W$ 
is a complex cone over $V$. We now compute the twisting of this 
fibration. The twisting is determined via the one-form 
$\diff\gamma+(n+1)(A+F(y)\diff\alpha)$ in the metric
(\ref{longmetric}). The integral of the corresponding curvature 
two-form is 
\bea\label{orbichern}
\frac{n+1}{2\pi}\int_W \diff F(y) \diff \alpha \, = \, \frac{p}{d(p-d)}~.
\eea
Fixing a point on $V$, and fixing $r>0$, we obtain a circle 
orbibundle over $W$. The right hand side of (\ref{orbichern}) 
is minus the Chern number\footnote{For an explanation, see 
section \ref{bigorbs}.} of this 
orbibundle, and corresponds to the canonical line orbibundle over 
$W=\mathbb{WCP}^1_{[d,p-d]}$. This is given by 
$K_W = \mathcal{O}(-p)$. One way to see this is via the K\"ahler quotient description of 
the weighted projective space together with its canonical line bundle 
over it. This is $\C^3//U(1)$ 
where the $U(1)$ action has weights $(d,p-d,-p)$. The weighted projective 
space itself is $z_3=0$, in standard complex coordinates on $\C^3$.

The above fibration structure immediately implies the earlier statements 
about the period of $\gamma$ over the singular points of $W$. In order 
to see this one needs to know some facts about orbibundles and 
orbifold fibrations. Suppose $W$ is an orbifold, with local orbifold 
charts $\{U_{\alpha} = \tilde{U}_{\alpha}/\Gamma_{\alpha}\}$, where 
$\tilde{U}_{\alpha}$ is an open set in $\R^N$ and $\Gamma_{\alpha}$ is 
a finite subgroup of $GL(N,\R)$. The data that defines an orbibundle 
over $W$ with structure group $G$ includes elements 
$h_{\alpha}\in \mathrm{Hom}(\Gamma_{\alpha},G)$ for each $\alpha$, subject to 
certain gluing conditions. In particular, if $F$ denotes a fibre 
over a smooth point of $W$, on which $G$ acts, 
then the fibre over a singular point with 
orbifold structure group $\Gamma_{\alpha}$ is 
$F/h_{\alpha}(\Gamma_{\alpha})$. Thus an orbibundle is generally 
not a fibration in the usual sense, since not all fibres are isomorphic. 

In the present situation it is particularly simple to work out 
the representations $h_{\alpha}$, since the orbibundle we require 
is the canonical line bundle $K_W$ over $W$. Since $W$ is a complex orbifold 
of dimension one, 
this is the holomorphic cotangent orbibundle. The orbifold structure groups 
are of the form $\Z_d,\Z_{p-d}\subset U(1)$, and then the maps $h_1: 
\Z_d\rightarrow U(1)$, $h_2: \Z_{p-d}\rightarrow U(1)$ are just the 
standard embeddings into $U(1)$. 
This implies that the metrics above are defined on 
\bea
K_{\mathbb{WCP}^1_{[d,p-d]}}\times_{U(1)} \C_V~.\eea
Here the $U(1)\subset\C^*$ action is the standard one on the canonical complex cone 
$\C_V$. The fibres over the poles of $W=\mathbb{WCP}^1_{[d,p-d]}$ 
are then $\C_V/\Z_d$ and $\C_V/\Z_{p-d}$, where the cyclic 
groups are embedded in $U(1)$ in the standard way. Here 
we have used the above maps $h_{\alpha}$, $\alpha=1,2$.
This completes the proof of Theorem \ref{flopper}.

Note that, in contrast to the previous section, we only obtain 
metrics for which $k=Id$ is divisible by $I$. Note also that the 
weighted projective space is a smooth $\mathbb{CP}^1$ if and only if 
$p=2$, $k=I$, which is a limiting case of the solutions considered 
here.


\section{Canonical resolutions}
\label{bigones}

In this section we turn our attention to complete 
Ricci-flat K\"ahler orbifold metrics, where the conical 
singularity gets replaced by a divisor $M$ with at worst
orbifold singularities. In section
 \ref{regularitysec} we addressed regularity of the metrics 
for $\pm x > \pm x_{\pm}$, and this 
 fixed uniquely the value of the parameter $\nu$ in terms of the pair 
of integers $p$ and $k$, in the range
 (\ref{pkrange}). The strategy here will be to show that one can choose 
appropriate values for the parameter 
 $\mu$ so that the metrics collapse smoothly, 
in an orbifold sense, to a divisor $M$ at $x=x_+$ or $x=x_-$, 
provided
$x_+>1$ and $x_-<y_1$. 
In fact for each $p$ and $k$ we shall find a family of values of 
$\mu$, indexed by an integer $r$ with $0<r<k/I$. $M$ is then a Fano 
orbifold of complex dimension $n+1$ which is a 
$\mathbb{WCP}^1_{[r,p-r]}$ fibration 
over $V$. The Ricci-flat K\"ahler metric is defined on the total 
space of the canonical 
line bundle over $M$. The induced metric on $M$ is
 K\"ahler, though in general the K\"ahler class is irrational. 
In order that $M$ be smooth ones requires $p=2$, $r=1$, and this leads to 
 Corollary \ref{ruled}.

\subsection{Partial resolutions III}\label{bigorbs}

Again, due to the symmetry in 
$x$ and $y$, the analysis of the collapse at $x=x_{\pm}$, for $y_1<y<y_2$, 
is identical to that of the collapse at $y=y_i$ for $\pm x>\pm x_{\pm}$. 
Thus for fixed $y\in (y_1,y_2)$ we deduce that the metric 
collapses smoothly at $x=x_{\pm}$ for all $\mu$. It thus remains to check 
the behaviour of the metric at $\{x=x_{\pm},y=y_1\}$ and 
$\{x=x_{\pm},y=y_2\}$.

We begin by writing the induced metric on $x=x_{\pm}$
\bea\label{inducedM}
\mp g_M &= &(1-x_{\pm})(1-y)g_V + \frac{y-x_{\pm}}{4Y(y)}\diff y^2 \nn\\
&+ &  \frac{x_{\pm}^2 Y(y)}{y-x_{\pm}}\left[\diff 
\alpha + \frac{x_{\pm}-1}{(n+1)x_{\pm}}[\diff\gamma\!+\!(n+1)A]\right]^2 \!\!.\eea
Near to a root $y=y_i$ we define
\bea
\mp R_i^2 = \frac{(y_i-x_{\pm})(y-y_i)}{Y'(y_i)}~,\eea
so that  near each  root we have the positive definite metric
\bea
g_M &  = & \pm(x_{\pm}-1)(1-y_i +\mathcal{O}(R_i^2))g_V +\diff R_i^2  \nn\\
& + & \left[\frac{(n+1)x_{\pm}y_iR_i}{y_i-x_{\pm}}\right]^2\left[
\diff\alpha+ \frac{x_{\pm}-1}{(n+1)x_{\pm}}[\diff\gamma+(n+1)A]\right]^2+{\cal O}(R_i^4)~.\eea
Let us define
\bea
\varphi_i \, =\,  \frac{(n+1)x_{\pm}y_i}{y_i-x_{\pm}}\left[\alpha+\frac{x_{\pm}-1}{(n+1)x_{\pm}}\gamma\right]~.
\label{beaver}
\eea
In order to allow for orbifold singularities, we impose the periodicities
\bea
\Delta \varphi_1  \,=\,  \frac{2\pi}{r}\qquad\quad \Delta \varphi_2  \,=\,  \frac{2\pi}{s}
\label{stinky}
\eea
for $r,s$ positive integers. This implies the necessary condition
\bea
-\frac{ry_1}{y_1-x_{\pm}} &= &\frac{sy_2}{y_2-x_{\pm}}
\label{pussycat}
\eea
where the minus sign ensures that both sides of the equation have the same sign.
This then  gives
\bea
x_{\pm} \, =\,  \frac{(r+s)y_1y_2}{ry_1+sy_2}~.
\label{panic}\eea
We shall return to this formula in a moment. 
A calculation using (\ref{panic}) and 
(\ref{bundletwists}) shows that
\bea
\varphi_1& =& \left(1+\frac{s}{r}\right)\frac{1}{p}\tilde{\alpha} + \frac{k(1+\tfrac{s}{r})-pI}{pI}\gamma\nn\\
-\varphi_2& = &\left(1+\frac{r}{s}\right)\frac{1}{p}\tilde{\alpha} + \frac{k(1+\tfrac{r}{s})-pI}{pI}\gamma~.
\label{winnythepooh}
\eea
Recall that $\tilde \alpha$ and $\gamma$ have period $2\pi$. 
In order to satisfy (\ref{stinky}) we must then 
require that $p=s+r$, which gives
\bea
\varphi_1 & =& \frac{1}{r}\tilde{\alpha} + \left(\frac{k}{rI}-1\right)\gamma\nn\\
-\varphi_2 &= &\frac{1}{s}\tilde{\alpha} + \left(\frac{k}{sI}-1\right)\gamma~.
\label{ribbon}
\eea
Let us now examine (\ref{panic}). Since the numerator is negative 
definite, for $\nu\in (0,\numax)$, this implies that
\bea
x_- &= &\frac{(r+s)y_1y_2}{ry_1+sy_2}~, \qquad \mathrm{for}\qquad r y_1+sy_2>0\nn\\
x_+ &=&\frac{(r+s)y_1y_2}{ry_1+sy_2}~, \qquad \mathrm{for}\qquad ry_1+sy_2<0~.
\eea
Note  that $ry_1+sy_2=0$ implies from (\ref{pussycat}) that $y_1=y_2$, which 
is impossible for $\nu\in (0,\numax)$. In particular we have
\bea
x_{-} - y_1 ~ = ~ \frac{ry_1(y_2-y_1)}{ry_1+sy_2} ~<~0
\eea
since the numerator is negative. 
By Lemma \ref{freak} there is therefore 
a unique $\mu$ such that $X(x)$ has $x_-$ as its smallest zero. 
On the other hand, using (\ref{bundletwists}) it is easy to compute
\bea\label{kylie}
\frac{y_1(x_+-1)}{y_1-x_+} = \frac{k-pI}{pI}+\frac{s}{r}\frac{k}{pI}~.\eea
Since $y_1/(y_1-x_+)$ is certainly positive, this implies that 
$x_+>1$ if and only if the right hand side of (\ref{kylie}) is 
positive. Using $r+s=p$, this easily becomes
\bea
\label{inbox}
x_+>1 \qquad \mathrm{iff} \qquad k-rI>0~.\eea
Thus when $ry_1+sy_2<0$ and $k-rI>0$, by Lemma \ref{ape} there is 
a unique $\mu<0$ such that $x_+>1$ is the largest zero of $X(x)$.

We now define
\bea
m \, = \, k-rI
\eea
and compute
\bea
ry_1+sy_2 \, =\,  \frac{1}{I}\left[(ky_1+(pI-k)y_2 + m(y_2-y_1)\right]~.\eea
Here we have substituted $s=p-r$. Using
\bea
y_1 &=& \left(\frac{k}{pI}-1\right)(y_2-y_1)+y_1y_2\nn\\
y_2 & = & \frac{k}{pI}(y_2-y_1)+y_1y_2~,\eea
which is a rewriting of (\ref{bundletwists}), we thus have
\bea\label{mboy}
ry_1+sy_2 \,=\, py_1y_2 + \frac{m}{I}(y_2-y_1)~.
\eea

Suppose that $m>0$. Then either $ry_1+sy_2>0$ and 
we are on the $x_-$ branch, with $x_-<y_1$; or else $ry_1+sy_2<0$ and 
by (\ref{inbox}) we are on the $x_+$ branch, with $x_+>1$. 
If $m <0$ then from (\ref{mboy}) $ry_1+sy_2<0$ and hence we are on the $x_+$ branch; 
but by (\ref{inbox}) $x_+<1$ and hence the metric cannot be regular. 
When $m=0$ we formally obtain $x_+=1$, which was the special case 
considered in the previous section.
We thus conclude that we obtain regular orbifold metrics if and only if $m>0$.

We have now shown that the metric $g_M$ extends to a smooth orbifold metric, 
for all $p,k,r$ positive integers with 
\bea
\frac{p}{2}<\frac{k}{I}<p, \qquad 0<r<\frac{k}{I}~.\eea
It 
remains simply to check the fibration structure of $M$ and thus describe its topology. Defining
\bea
\varphi \,=\, \frac{(n+1)ky_1}{I(y_1-1)}\left[\alpha+\frac{x_{\pm}-1}{(n+1)x_{\pm}}\gamma\right] \,=\,
 \tilde{\alpha} + \frac{k-rI}{I}\gamma\eea
we see that the one-form in the second line of the metric (\ref{inducedM}) is 
proportional to
\bea
\diff \varphi + \frac{k-rI}{I}(n+1)A~.\eea
Since $\varphi$ has canonical period $2\pi$, this is a global angular 
form on the associated circle bundle to $K_V^{m/I}$, where
recall $m=k-rI$. Thus $M$ may be described as follows. One takes the 
weighted projective space $\mathbb{WCP}^{1}_{[r,p-r]}$ and fibres this over
$V$. The transition functions are precisely those for $K_V^{m/I}$, using the 
standard effective $U(1)$ action on $\mathbb{WCP}^1_{[r,p-r]}$. 
Thus $M$ may be written 
\bea
M \, =\,  K_V^{m/I}\times_{U(1)}\mathbb{WCP}^1_{[r,p-r]}~.
\eea
The Ricci-flat K\"ahler metric is defined on the total space of an orbifold line bundle 
over $M$, which is necessarily the canonical line orbibundle.

Notice that $M$ is singular precisely along the two divisors 
$D_1$, $D_2$, located at $y=y_1$, $y=y_2$, respectively. 
$D_1$ has normal fibre $\C/\Z_r$, and $D_2$ has normal fibre 
$\C/\Z_{p-r}$. The normal bundles are $K_V^{-m/I}$, $K_V^{m/I}$, 
respectively. Due to the fact that the only orbifold singularities 
are in complex codimension one, $M$ is in fact completely smooth as a
manifold, and as an algebraic variety\footnote{Note that $V$ is 
a smooth Fano manifold, and hence is projective.}. In either case,
$M$ is a $\mathbb{CP}^1$ fibration over $V$. One must then be 
extremely careful when making statements such as ``M is Fano'': 
the anti-canonical line bundle and anti-canonical orbifold line bundle 
are different objects. 

Let $\pi:M\rightarrow V$ denote the projection. Then the 
canonical line bundle is
\bea\label{canonicalclass}
K_M \,=\, \pi^*K_V -2D_1 -m \pi^*(K_V/I)~.\eea
Recall here that the divisor $D_1$ at $y=y_1$ 
has normal bundle $K_V^{-m/I}$. 
Note in (\ref{canonicalclass}) 
we have switched to an additive notation, rather than 
the multiplicative notation we have been using so far throughout 
the paper; this is simply so that the equations are easier to read. 
On the other hand, the orbifold canonical line bundle is
\bea\label{orbi}
K_M^{\mathrm{orb}} \,=\, K_M + \left(1-\frac{1}{r}\right)D_1 + \left(1-\frac{1}{p-r}\right)D_2~.\eea
This may be argued simply by the following 
computation, taken largely from \cite{kollarpaper}. 
Let $U$ be an open set in $M$ 
containing some part of a divisor $D$ with normal fibre $\C/\Z_r$. 
We suppose that $\tilde{U}\subset \C^n$ is the local covering chart, 
and that the preimage of the divisor $D$ is given locally in $\tilde{U}$ 
by $x_1=0$. This is called the ramification divisor; we  
denote this divisor in $\tilde{U}$ by $R$. We also complete $x_1$ 
to a set of local complex coordinates on
$\tilde{U}$, $(x_1,\ldots,x_n)$. The orbifold structure group is 
$\Gamma=\Z_r$, and the map $\phi:\tilde{U}\rightarrow U$ near $D$ looks like
\bea
\phi: (x_1,x_2,\ldots,x_n)\rightarrow (z_1=x_1^{r},z_2=x_2,\ldots,z_n=x_n)\eea
where $(z_1,\ldots,z_n)$ are complex coordinates on $U$, which is 
also biholomorphic to an open set in $\C^n$. In particular, we may compute
\bea\label{covering}
\phi^*(\diff z_1\wedge\cdots\wedge\diff z_n) \,=\, rx_1^{r-1}\diff x_1 
\wedge\cdots\wedge\diff x_n~.\eea
Now, the orbifold line bundle $K_M^{\mathrm{orb}}$ is 
defined as the canonical line bundle of $\tilde{U}$ in each covering 
chart $\tilde{U}$, {\it i.e.} as the top exterior power of the 
holomorphic cotangent bundle. These naturally glue together on 
the orbifold $M$ to give an orbifold line bundle over $M$. 
However, we see from (\ref{covering}) that
\bea
K_{\tilde{U}} \,=\, \pi^*K_{U}\otimes [(r-1)R]~.\eea
Since $\pi^*D = rR$, this gives the general formula
\bea\label{formal}
K_{M}^{\mathrm{orb}} \,=\, K_M + \sum \left(1-\frac{1}{r_i}\right)D_i\eea
where $D_i$ is a so-called branched divisor, with multiplicity $r_i$.  
This rather formal expression may be understood more 
concretely as follows. $M$ is a complex manifold with divisors 
$K_M$ and $D_i$ defining complex line bundles over $M$. For each line 
bundle, by picking a 
connection we obtain a curvature two-form whose cohomology class
lies in the image of $H^2(M;\Z)$ in $H^2(M;\R)$. The corresponding cohomology 
class of the right hand side of (\ref{formal}) is thus in $H^2(M;\mathbb{Q})$. 
This in fact represents the cohomology class of the curvature of a 
connection on the \emph{orbifold} line bundle $K_M^{\mathrm{orb}}$. 

Returning to (\ref{orbi}), we obtain
\bea\label{neworbs}
K_M^{\mathrm{orb}} &= &-\frac{p}{r(p-r)} D_1 + \left(I-\frac{m}{p-r}\right)\pi^*(K_V/I)\nn\\
&=& -\frac{p}{r(p-r)} D_2 + \left(I+\frac{m}{r}\right)\pi^*(K_V/I)
\eea
where note that
\bea
D_1-D_2\, =\, -m\pi^*(K_V/I)~.
\eea
Note in (\ref{neworbs}) the first Chern class of the weighted projective space 
$\mathbb{WCP}^1_{[r,p-r]}$ appearing. Indeed, by the above comments 
the integral of this orbifold first Chern class is
\bea
\int_{\mathbb{WCP}^1_{[r,p-r]}} c_1^{\mathrm{orb}}\left(\mathbb{WCP}^1_{[r,p-r]}\right) \, =\, 
2 - \left(1-\frac{1}{r}\right) 
- \left(1-\frac{1}{p-r}\right) = \frac{p}{r(p-r)}~,\eea
a formula we encountered earlier in equation (\ref{orbichern}).

Let $\Sigma\subset V$ be a holomorphic curve in $V$. We may map 
$\Sigma$ into $M$ via the sections $s_i:V\rightarrow M$ at $y=y_i$. 
Using (\ref{neworbs}) we then compute
\bea
\langle c_1(K_M^{\mathrm{orb}}),s_1(\Sigma)\rangle & = & \left[\frac{p}{r(p-r)}m + \left(I-\frac{m}{p-r}\right)\right]\langle c_1(K_V^{1/I}),\Sigma\rangle\nn\\
& =& \frac{k}{r}\langle c_1(K_V^{1/I}),\Sigma\rangle~,
\eea
\bea
\langle c_1(K_M^{\mathrm{orb}}),s_2(\Sigma)\rangle & = & \left[-\frac{p}{r(p-r)}m + \left(I+\frac{m}{r}\right)\right]\langle c_1(K_V^{1/I}),\Sigma\rangle\nn\\
& =& \frac{pI-k}{p-r}\langle c_1(K_V^{1/I}),\Sigma\rangle~.\eea
Here we have used that {\it e.g} $\langle D_1, s_1(\Sigma)\rangle = -m 
\langle c_1(K_V^{1/I}),\Sigma\rangle$, since $y=y_1$ has normal bundle 
$K_V^{-m/I}$. Since $V$ is Fano, $\langle c_1(K_V^{1/I}),\Sigma\rangle<0$, and hence $M$ 
Fano implies that 
\bea
k-pI<0~.
\eea
Of course, this condition is indeed satisfied by the explicit 
metrics we have constructed. This completes the proof 
of Theorem \ref{bigorbifolds}.

\subsection{Smooth resolutions: $p=2$}

Setting $p=2$, $r=1$ in the last subsection gives a family of 
smooth complete Ricci-flat K\"ahler metrics for each choice of $(V,g_V)$, 
leading to Corollary \ref{ruled}. 
These are all defined on the canonical line bundle over $M$, where
$M=\mathbb{P}_V(\mathcal{O}\oplus K^{m/I})$, $m=k-I$, and 
$0<m<I$. 

For example, we may take $V=\mathbb{CP}^n$ with its standard K\"ahler-Einstein
metric. In this case  $0<m<n+1$ and the 
metrics are defined on the total space of the canonical 
line bundle over $\mathbb{P}_{\mathbb{CP}^n}(\mathcal{O}(0)\oplus\mathcal{O}(-m))$. Note that $n=m=1$ is precisely the metric found in reference 
\cite{japanese}.

On the other hand, taking $V$ to be a product of 
complex projective spaces, as in (\ref{product}), 
reproduces the metrics discussed in reference \cite{lupope}. 
In the latter reference the authors considered each product separately; 
this was necessary, given the method they use to analyse regularity 
of the metric. The number of smooth metrics found is $I-1=\mathrm{hcf}\{d_a\}-1$; one can easily verify that 
this number agrees with the number of smooth resolutions found in 
the various cases considered in \cite{lupope}.


\subsection*{Acknowledgments}
\noindent We thank S.-T. Yau for discussions on 
related topics. J. F. S. would also
like to thank the mathematics department at the University of California, 
Los Angeles, for hospitality. He
is supported by NSF grants DMS-0244464, DMS-0074329 and DMS-9803347.
D. M. would like to thank the physics and mathematics departments of 
Harvard University for hospitality during completion of this work. 
He acknowledges support from NSF grant PHY-0503584.

\appendix

\section{Limits}

In this section we briefly analyse various special limits of
 the metrics (\ref{metric}). Recall that 
these depend on four real parameters:
 $c,\beta,\nu,\mu$. 
First, note that setting $c=0$ in the functions (\ref{functions}) implies
that  the $g_{\tau \tau}$ component
of the metric asymptotes to a constant ($\mp\beta$) as $x\to \pm\infty$.
 Therefore the metric is not asymptotically conical. 
When $c\neq0$ we may then set $c=1$
by a diffeomorphism and rescaling of the metric, 
as we have assumed throughout the paper.

If $\beta$ is different from zero, 
it may also be scaled out as an overall coefficient of the metric 
(\ref{metric}), where $\beta=1$.
However, we may also consider asymmetric scalings of the variables 
$x$ and $y$, before 
letting $\beta\to 0$, with the result depending on which variable 
goes to zero faster. There are then two cases 
to consider. Thus, let us first make the substitution 
$y  \to  \beta y $. The resulting metric reads
\bea
g & = & \frac{\beta y-x}{4X_\beta (x)}\diff x^2 + \frac{\beta y-x}{4Y(y)}\diff y^2 + \frac{X_\beta(x)}
{\beta y-x}
\left[\diff\tau + (1-y)(\diff\psi + A)\right]^2 \nonumber\\
&& + \frac{Y(y)}{\beta y-x}\left[\beta \diff\tau + (\beta-x)(\diff\psi + A)\right]^2 +
(\beta-x)(1-y)g_V~, \label{betametric}
\eea
where the parameter $\nu$ in $Y(y)$ has been redefined, and we have introduced the notation 
\bea
X_\beta (x) & = & \beta(x-\beta)+\frac{n+1}{n+2}(x-\beta)^2+\frac{2\mu}{(x-\beta)^n}
\eea
to emphasize that $X_\beta(x)$ depends $\beta$, as opposed to $Y(y)$. Setting $\beta=0$
in (\ref{betametric}) and introducing the change of variable
\bea
x & = & \pm \frac{n+1}{n+2}r^2
\eea
as in section \ref{asymptotics}, we obtain the positive definite metric
\bea
g & = & \frac{1}{H(r)}\diff r^2 + H (r)r^2 \left(\frac{n+1}{n+2}\diff \tau +\sigma  \right)^2+ r^2 g_T~,
\eea
where 
\bea
H(r)& = & 1 + 2\mu\left(\frac{n+1}{n+2}\right)^{n+3}\frac{(-1)^n}{r^{2n+4}}~.
\label{PPrecovered}
\eea
This is precisely the Calabi ansatz of \cite{BB,PP}. Note from the metric (\ref{betametric}) 
that $\beta$ plays the role of a resolution parameter
in all cases considered in the paper.  The parameter $\nu$ is 
fixed by regularity at 
$x>x_+$ or $x<x_-$, as discussed in section \ref{regularitysec}, 
and one is left with two paramaters $\mu$ and $\beta$. 
This is analogous to the two-parameter family of Ricci-flat metrics on the 
canonical line bundle
over $\mathbb{CP}^1\times \mathbb{CP}^1$ found in \cite{Pando-Tseytlin}, 
which has regular 
asymptotic boundary metric $T^{1,1}/\Z_2$. However, in the 
present case, regularity (in the orbifold sense)
of the metrics at $x=x_{\pm}$
 imposes a relation. This leaves $\beta$ as the only free parameter, 
measuring the size of the blown-up cycles. 
It would be interesting to investigate
generalisations of the ansatz (\ref{metric}), allowing for more than one 
 resolution parameter.

Finally, it is straightforward to repeat the previous analysis setting instead $x\to \beta x$, 
by exchanging the 
roles of $x$ and $y$. However, in this case if we set $\beta = 0$, the 
metric becomes 
\bea
g & = & \frac{y}{4X(x)}\diff x^2 + \frac{y}{4Y_0(y)}\diff y^2 + X(x)y(\diff\psi + A)^2 \nn\\
& + & \frac{Y_0(y)}{y}\left[\diff\tau + (1-x)(\diff\psi + A)\right]^2 -y(1-x)g_V~, 
\eea
where $Y_0(y)$ denotes the function $Y(y)$ in (\ref{functions}), evaluated at $\beta=0$.
We see  the metric is again not asymptotically conical.

\end{document}